\documentclass[]{interact}
\usepackage{epstopdf}

\usepackage[numbers,sort&compress]{natbib}
\bibpunct[, ]{[}{]}{,}{n}{,}{,}
\makeatletter
\def\NAT@def@citea{\def\@citea{\NAT@separator}}
\makeatother

\usepackage[latin1]{inputenc}
\usepackage[english]{babel}

\usepackage[makeroom]{cancel}
\usepackage{xfrac}
\usepackage{ulem}
\usepackage{color}%
\usepackage{amsmath}%

\usepackage[hidelinks]{hyperref}
\usepackage{tensor} 

\newcommand{\dejak}[3]{\tensor[^{#1}]{\overset{\curvearrowright}{#2}}{^{#3}}}

\usepackage{amsfonts}%
\usepackage{amssymb}
\usepackage{mathabx}
\usepackage{comment}
\usepackage{graphicx}
\usepackage{subfigure}
\usepackage{t1enc,epsfig,multicol}
\DeclareMathOperator{\dist}{dist}

\DeclareMathOperator{\Diag}{Diag}

\DeclareMathOperator{\diag}{diag}

\DeclareMathOperator{\conv}{conv}

\newcommand{\wjk}{\widetilde{w}^{jk}}

\def\R{\mathbb{R}}   
\def\C{\mathbb{C}}   
\def\V{\mathbb{V}} 
\def\W{\mathbb{W}} 
\def\DD{\mathcal{D}}
\def\DV{\mathcal{D}_{\V}}
\def\DS{\mathcal{D}_{S}}
	\def\mV{m_{\mathbb{V}}}   %
	\def\mS{m_{S}}   %
\def\tr{\text{tr}}
 
\usepackage[numbers,sort&compress]{natbib}
\bibpunct[, ]{[}{]}{,}{n}{,}{,}
\makeatletter
\def\NAT@def@citea{\def\@citea{\NAT@separator}}
\makeatother

\theoremstyle{plain}
\newtheorem{theorem}{Theorem}[section]
\newtheorem{lemma}[theorem]{Lemma}
\newtheorem{corollary}[theorem]{Corollary}
\newtheorem{proposition}[theorem]{Proposition}

\theoremstyle{definition}
\newtheorem{definition}[theorem]{Definition}
\newtheorem{example}[theorem]{Example}

\theoremstyle{remark}
\newtheorem{remark}{Remark}
\newtheorem{notation}{Notation}
\begin{document}
%


\title{Moment of a subspace and joint numerical range}

\author{
	\name{Abel Klobouk\textsuperscript{a} and Alejandro Varela\textsuperscript{b}\thanks{CONTACT A. Varela.
		Email: avarela@ungs.edu.ar}}
	\affil{\textsuperscript{a}Departamento de Ciencias B{\'a}sicas, Universidad Nacional Luj{\'a}n, Luj{\'a}n, Argentina; \textsuperscript{b}Instituto de Ciencias, Universidad Nacional de General Sarmiento, Los Polvorines, Argentina} }

\maketitle

\begin{abstract}
	For a given complex finite dimensional subspace $S$ of $\C^n$ and a fixed basis, we study the compact and convex subset of $\left(\R_{\geq 0}\right)^n$ that we call the moment of $S$
	\begin{equation*}
	\begin{split}
		m_S=&\text{convex } \text{hull}\left(\{|s|^2\in\R^n_{\geq 0}:s\in S\wedge \|s\|=1\}\right)\\
		    \simeq & \{{\Diag}(Y) \in M_n^h(\C):Y\geq 0, \tr(Y)=1, P_S Y P_S=Y\}
	\end{split}
	\end{equation*}
	where $|s|^2=(|s_1|^2,|s_2|^2,\dots,|s_n|^2)$. This set is relevant in the determination of minimal hermitian matrices ($M\in M^h_n$ such that $\|M+D\|\leq D$ for every diagonal $D$ and $\|\ \|$ the spectral norm). We describe extremal points and curves of $m_S$ in terms of principal vectors that minimize the angle between $S$ and the coordinate axes. We also relate $m_S$ to the joint numerical range $W$ of $n$ rank one $n\times n$ matrices constructed with the orthogonal projection $P_S$ and the fixed basis used. This connection provides a new approach to the description of $m_S$ and to minimal matrices. As a consequence the intersection of two of these joint numerical ranges allow the construction or detection of a minimal matrix, a fact that is easier to corroborate than the equivalent condition for moments. It is also proved that $m_S$ is a semi-algebraic set equal to the intersection of the mentioned $W$ with a hyperplane and whose generated positive cone coincides with that of $W$.
\end{abstract}


\begin{amscode}{65F35, 15B57, 15A60, 52A20.} \end{amscode}

\section{Introduction}
Given a subspace $S$ of $\C^n$, our main interest in the study of the set
\[
m_S=\conv\left(\{|s|^2\in\R^n_{\geq 0}:s\in S\wedge \|s\|=1\}\right)
\]
where $|s|^2=\left(|s_1|^2, |s_2|^2,\dots,|s_n|^2\right)$ with $s_j$, $j=1,\dots, n$ are the coordinates of $s$ in a fixed basis $E=\{e_1, e_2, \dots, e_n\}$ (see Proposition \ref{prop: equivalencias de momento} for other equivalent definitions of $m_S$).
Its importance lies on its fundamental relation with what we call minimal hermitian matrices $M\in M_n^(\C)$, that satisfy $\|M\|\leq \|M+D\|$, with $D$ a diagonal matrix and $\| M \|=\sup_{\|x\|=1}\|M x\|$. Note that the consideration of diagonal matrices $D$ requires the choice of a fixed orthonormal basis $E$ of $\C^n$. The relation between minimal matrices and $m_S$ is the following
\begin{equation}
\label{eq: relac momento y minimales}
M \text{ is minimal } \Leftrightarrow M=\|M\| \left(P_{S_{\|M\|}} -P_{S_{-\|M\|}} \right)+R
\ \text{ and } \ 
m_{S_{\|M\|}}\cap m_{S_{-\|M\|}}\neq \emptyset 
\end{equation}
where $P_S$ denotes de orthogonal projection onto $S$ and $\|R\|<\|M\|$, Im$(R)\perp S_{\|M\|}, S_{-\|M\|}$, with $S_{\lambda}$ denoting the eigenspaces of $M$ corresponding to the eigenvalue $\lambda$ (see the more detailed Remark \ref{rem: intersecc momentos entonces construcc mat minimal}). These matrices, in turn, allow the concrete description of some metric geodesics in flag manifolds as studied in \cite{ammrv, duranmatarecht}.

We describe many extremal points of $m_S$ (including curves of them) that can be easily described in terms of some principal vectors of $S$ or the matrix of the orthogonal projection $P_S$ in the $E$ basis. The projections of the mentioned curves of $m_S$ in certain 2 dimensional coordinate planes are parts of ellipses centred at the origin that can be easily obtained.

We also study a very close relationship between the moment $m_S$ and a joint numerical range 
of some chosen hermitian matrices of rank one $P_S E_i P_S=(P_S e_i)\cdot(P_S e_i)^*=(P_S e_i)\otimes(P_S e_i)$, for $i=1,\dots, n$ (see Theorem \ref{teo mV es JNR y suma 1}). This allows a translation of the condition of being a minimal matrix \eqref{eq: relac momento y minimales} to a property of intersection of the two corresponding joint numerical ranges (see Theorem \ref{teo: intersecc mS no vacia igual a intersec JNR no nula}). Recent results using algebraic geometry allow a description of a subset $T^\sim$ whose convex hull gives these joint numerical ranges (see \cite{plaumann-sinn-weis}). This is a generalization of a Theorem of Kippenhahn where it is proved that the convex hull of certain algebraic curve equals the classic numerical range of a matrix $W(A)$. 

In what follows we list the main contents of each section of this work. In Section \ref{sec: Preliminares} we state most of the notation used. Section \ref{moment of subspace} includes previous results and different characterizations related to $m_S$. We also introduce there what could be considered the centre of $m_S$ (see \eqref{def centroide}) and some of its properties.

In Section \ref{Generic subspaces and its principal standard vectors} we introduce the principal standard vectors $\{v^j\}_{j=1}^n$ of a subspace $S$ related to a fixed basis $E$. These are the ones that minimize the angle between $S$ and the coordinate axis of a fixed basis $E=\{e_1, e_2,\dots, e_n\}$. We describe some of its properties, particularly the ones related to the elements $|v^j|^2\in m_S$, $j=1,\dots,n$.

Section \ref{curves of extreme points} is dedicated to the presentation and description of certain particular curves of extreme points of $m_S$ that join elements $|v^j|^2 \in m_S$ of what we call the principal standard vectors $\{v^j\}_{j=1}^n$.

In Section \ref{secc JNR} we show the close relation between $m_S$ and the joint numerical range 
\begin{equation*}
\begin{split}
W_{S,E}&=W(P_S E_1 P_S, \dots, P_S E_n P_S)\\
&=\left\{\left(\tr(P_S E_1 P_S\rho, \dots,  \tr(P_S E_1 P_S\rho\right)\in \R^n_\geq 0: \rho\geq 0, \tr(\rho)=1\right\},
\end{split}
\end{equation*}
where $E_i=e_i\otimes e_i$ is the orthogonal projection onto the subspace generated by $e_i\in E$, $i=1,\dots n$. Here $E$ is the fixed standard basis we are considering. This allows restating many results of $m_S$ and minimal matrices in terms of these joint numerical ranges and obtain some general properties about $m_S$.
Among them it is proved that $m_S$ is a semi-algebraic set (Remark \ref{rem propiedades mS a partir de W} 
 c)), it is a subset of $W_{S,E}$ included in a hyperplane that generates the same cone than $W_{S,E}$ (Theorem \ref{teo mV es JNR y suma 1} and Proposition \ref{prop cono igual cono}) and the known algorithms used to generate $W_{S,E}$ can be used to approximate $m_S$. It is also shown that the map $P_S\mapsto m_S$ satisfies $\dist_H(m_S,m_V) \leq  (2 \sqrt{n}+1) \|P_S-P_V\|$ (Proposition \ref{prop cono igual cono}) where $\dist_H$ is the Hausdorff distance. Moreover,  
  for a hermitian matrix $M=\|M\| (P_{\text{Eig}_{\|M\|}}-P_{\text{Eig}_{-\|M\|}})+R$ (with ${\text{Eig}_{\|M\|}}$, ${\text{Eig}_{-\|M\|}}$ the eigenspaces of $\pm\|M\|$, $R\in M_n^h(\C)$,  $R({Eig_{\pm\|M\|}})=0$), holds that  $W_{{\text{Eig}_{\|M\|}},E}\cap W_{{\text{Eig}_{-\|M\|}},E}\neq \{0\}$ if and only if $M$ is minimal.
 
\section{Preliminaries and notation}\label{sec: Preliminares}

Let $\C^n$ be the finite $n$-dimensional Hilbert space of vectors of complex numbers with its usual inner product denoted with $\langle\ ,\ \rangle$ and norm $\|\ \|$.\
We denote with $E=\{e_1,e_2,\dots,e_n\}$ a fixed orthonormal ordered basis of $\C^n$ which we will call the  standard basis. The vectors $e_i\in E$, $i=1,\dots,n$ will be called the standard vectors.  Our adoption of a fixed basis is required since we are going to work with diagonal matrices that may not remain as such when a different basis is used. 

Given a vector $x\in\C^n$ we will denote with $x_1$, $x_2$, $\dots$,$x_n$ its standard coordinates (in the $E$ basis) with the only exception of $\{e_1, e_2,\dots,e_n\}$ where the subscripts denote the different standard vectors.

We will denote with $|x|^2=(|x_1|^2,\dots,|x_n|^2)\in\R^n_{\geq 0}$ the vector of the squared modulus of its $E$ coordinates $x_i=\langle x, e_i\rangle$, for $i=1,\dots,n$. The vector $|x|^2$ has been described in previous works with $x\circ \bar{x}$ using the entry-wise Schur (or Hadamard) product $\circ$ where $\bar{x}$ denotes the vector whose coordinates are those of $x$ conjugated. We state here some of these notation for further references.
\begin{notation}\label{notaciones} We will use the following notations:
	\begin{enumerate}
		\item[a)] $E=\{e_1, e_2, \dots, e_n\}$ will denote a fixed orthonormal ordered basis in $\C^n$, considered the standard basis,
		\item[b)] given $x\in\C^n$, the $x_j$ will be the $j^{\text{th}}$ coordinate of $x$ in the standard basis, that is, $x_j = \langle x , e_j \rangle$,
		\item[c)] \label{notacion |x| al cuad} given a vector $x\in\C^n$ we denote with $|x|$ the vector of the modules of its coordinates $|x|=(|x_1|,\dots,|x_n|)$ and similarly $|x|^2=(|x_1|^2,\dots,|x_n|^2)$.
	\end{enumerate}
\end{notation}

With $M_n(\C)$ we denote the set of $n\times n$ matrices with coefficients in $\C$, and if $A\in M_n(\C)$ then $A^*$ is its adjoint matrix (conjugate transpose of $A$). We write $A\in M_n^h(\C)$ if $A$ is a matrix such that $A=A^*$, and we will call it hermitian. Given a square matrix $A$, then $\Diag(A)$ is the diagonal square matrix with the same diagonal than $A$ written using a fixed basis of $\C^n$. Sometimes we will identify the matrix $\Diag(A)$ with the vector whose entries are obtained from the corresponding diagonal.

The rank one orthogonal projections onto the subspace generated by an $x\in \C^n$, with $\|x\|=1$ will be denoted with $x\otimes x$ ($(x\otimes x)(v)=\langle v,x\rangle x$, for all $v\in\C^n$). Its corresponding matrix in the $E$ basis can be written by $x\cdot x^*$ where $\cdot$ is the matrix product and $x$ is identified with the $n\times 1$ column matrix of coordinates of $x$ in the $E$ basis ($x_{i,1}=\langle x, e_i\rangle$).

If $T\subset \C^n$ or $\R^n$, we will write $\conv(T)$ to denote the convex closure or convex hull of $T$. Equivalently the set of convex combinations $t\, a+(1-t) b$ for $t\in[0,1]$, for all $a,b\in T$.

\bigskip

\section{Moment of a subspace}\label{moment of subspace}

The term ``moment'' in the following definition of the set $m_S$ of a subspace $S$ of dim$(S)=r$ is motivated by its relation to certain moment map defined in the symplectic manifold $(\C^n)^r$ (see Section 4 in \cite{soportes} for details).

\begin{definition}\label{def: momento de un subespacio} 
	Let $S$ a subspace of $\C^n$, $\{0\}\subsetneq S\subsetneq \C^n$ with $\dim(S)=r$.
The moment of $S$ related to a fixed orthonormal ordered basis $E=\{e_1,e_2,\dots,e_n\}$ of $\C^n$ or just the \textbf{moment set of $S$} is the subset of $\left(\R_{\geq 0}\right)^n$ defined by
\begin{equation}
m_{S,E}
= \conv\big( \{|v|^2\in\R^n_{\geq 0}: v\in S \wedge \|v\|=1\} \big).
\end{equation}	
If the basis $E$ is fixed or clear from the context we will use the shorter notation $m_S$ instead of $m_{S,E}$.
\end{definition}

\begin{remark} The set $m_S$ is a compact and convex set of $\R^n$, but note that if the convex closure ($\conv$) is not used in \eqref{def: momento de un subespacio} then $m_S$ would not be necessarily convex (see Remark 4 in \cite{alrv}). 
	
There are several equivalent statements that we can consider to define the set $m_S$ some of which do not require to take the convex hull. See for example Proposition \ref{prop: equivalencias de momento}.

\end{remark}
\begin{remark}	\label{rem: intersecc momentos entonces construcc mat minimal}
	Our motivation to study the set $m_S$ relies on the following property. Let $\V$ and $\W$ be two non trivial orthogonal subspaces of $\C^n$, and a hermitian matrix $R\in M^h_n(\C)$ such that $R P_\V=R P_\W=0$ and $\|R\|<\lambda$. Then,	if $m_\V\cap m_\W \neq  \emptyset$, follows that the matrix  $M=\lambda(P_\V-P_\W) +R$ is a minimal hermitian matrix in the sense that 
\[
	\|M\|=\displaystyle{\dist\left(M, {\Diag}_n(\R)\right)}=\displaystyle{\dist\left(M, {\Diag}_n(\C)\right)}
\]
	for ${\Diag}_n(X)$ the diagonal matrices with coefficients in $X$. Moreover, every minimal matrix can be written in this way (see Corollary 3 of \cite{alrv} and Theorem 3 of \cite{soportes} to obtain a detailed study and other equivalent statements).
\end{remark}

\begin{remark}
We denote the positive cone generated by a set $X\subset\R^n$ as
	\begin{equation}\label{eq: def cono conj}
	\text{cone}\left(X\right)=\{t\, x\in\R^n: t\in\R_{\geq 0} \wedge x\in X\}.
	\end{equation}
	Note that we can consider the positive cone generated by $m_S$ in $\R^n_{\geq 0}$ instead of just $m_S$ in order to obtain the sufficient condition mentioned in the previous Remark \ref{rem: intersecc momentos entonces construcc mat minimal}. Nevertheless, we choose the present Definition \ref{def: momento de un subespacio} of $m_S$ in order to work with bounded sets. 
\end{remark}
%

The following are equivalent characterizations of the moment $m_S$. The statements \eqref{defequiv mS 1 en prop}, \eqref{defequiv mS 2 en prop},  \eqref{defequiv mS 3 en prop},  can be obtained from characterizations made in \cite{alrv}, \eqref{eq relac momento y trEiY} is Lemma \ref{lema: momento igual a trazas EiY} and \eqref{eq relac momento y JNR} is Theorem \ref{teo mV es JNR y suma 1} of section \ref{secc JNR}.
\begin{proposition}
	\label{prop: equivalencias de momento}
	The following are equivalent conditions to define $m_S$, the moment of $S$ related to a basis $E=\{e_1,\dots,e_n\}$, with dim$(S)=r$.
\begin{enumerate}
	\item \label{defequiv mS 1 en prop} Let $\Diag(X)\in M_n(\C)$ be the diagonal matrix with the same diagonal than $X \in M_n(\C)$ when its written using the basis $E$, then
		\[
		m_S  \simeq \{{\Diag}(Y) \in M_n^h(\C):Y\geq 0, \tr(Y)=1, \text{Im}(Y) \subset S\},
		\] 
		where we denote with $\simeq$ the usual identification of vectors of $\C^n$ with diagonal matrices of $\C^{n\times n}$. 
	\item \label{defequiv mS 2 en prop}
		\[
		m_S= {\conv}\left\{|v|^2: v \in S \wedge \|v\|=1 \right\}
		\]
		(see the meaning of $|v|^2$ in Notation \ref{notacion |x| al cuad})
	\item \label{defequiv mS 3 en prop}
		\[
		m_S= \bigcup_{\{s^i \}_{i=1}^r  \text{o.n. set in } S} \  {\conv}\left( \{|s^i|^2 \}_{i=1}^r\right).
		\]
	\item \label{eq relac momento y trEiY} 
		\[
		m_S=\{\left(\tr(E_1 Y), \tr(E_2 Y),\dots,\tr(E_n Y)\right)\in\R^n_{\geq 0}: Y \in \DD_S \}
		\] 
		where $E_i=e_i\cdot ( e_i )^t$ and $\DD_S= \{Y\in M_n^h(\C):Y\geq 0, \tr(Y)=1, P_S Y=Y\}$.
	\item \label{eq relac momento y JNR} 
		\[
		m_S=W(P_S E_1 P_S,\dots,P_S E_n P_S)\cap \left\{x\in \R^n_{\geq 0}:\sum_{i=1}^n x_i=1\right\},
		\] 
		where $P_S$ is the orthogonal projection onto $S$, $E_i=e_i\cdot(e_i)^*$ and $W$ is the joint numerical range defined in \ref{def JNR}.
\end{enumerate}
\end{proposition}
\begin{remark}
	Note that in particular the condition \eqref{defequiv mS 3 en prop} implies that the Carath{\'e}odory number of $\left\{|v|^2: v \in S \wedge \|v\|=1 \right\}$ is less or equal than $r=\dim(S)$ (see \eqref{defequiv mS 2 en prop}) .
\end{remark}
\begin{example} \label{ejemplo 1} Let $U_d\in M_n(\C)$ be a unitary such that its matrix with respect to the $E$ basis is diagonal, $S$ a subspace of $\C^n$ and $U_d(S)=\{U_d(s): s\in S \}$, the subspace that is the image of $S$ under $U_d$. Then it is trivial that $m_{S,E}=m_{U_d(S),E}$.
Nevertheless, if $\V, \W$ are subspaces of $\C^n$ such that $m_\V=m_\W$ it is not true that there exists a diagonal unitary $U_d$ (relative to the base $E$) such that $U_d(\V)=\W$. For instance, if we consider the subspaces $\V=\text{span}\{(1,1,0),(0,1,1)\}$ and $\W=\text{span}\{(-1,e^{{i \pi}/{4}},0),(0,e^{i \pi/3},e^{i \pi/6})\}$ it can be proved that $m_\V=m_\W$ but there is not any diagonal unitary $U_d$ such that $U_d(\V)=\W$. Their corresponding orthogonal projections written in the $E$ basis are 
\[
P_\V=
\begin{pmatrix}
 {2}/{3} &  {1}/{3} &-{1}/{3} \\
{1}/{3} & {2}/{3} & {1}/{3} \\
-{1}/{3} &  {1}/{3} & {2}/{3} \\
\end{pmatrix}\ \text{ and }\ 
P_\W=\begin{pmatrix}
\frac{2}{3} & -\frac{1}{3} e^{-\frac{i
		\pi }{4}} & \frac{1}{3} e^{-\frac{i
		\pi }{12}} \\
-\frac{1}{3} e^{\frac{i \pi }{4}} &
\frac{2}{3} & \frac{1}{3} e^{\frac{i
		\pi }{6}} \\
\frac{1}{3} e^{\frac{i \pi }{12}} &
\frac{1}{3} e^{-\frac{i \pi }{6}} &
\frac{2}{3}
\end{pmatrix}.
\]
\end{example}


\begin{remark}\label{rem igual diag para proy bon S}
	Let $S \subset \C^n$ be a subspace of dimension $r$ and $E$
	a fixed basis of $\C^n$.
		If $\{s^1,s^2,\dots,s^r\}$ and  $=\{w^1,w^2,\dots,w^r\}$ are two orthonormal basis of $S$, then 
		\begin{equation}\label{igualdad_de_centroides_expandidos}
		\sum_{i=1}^r |s^i|^2 = \sum_{i=1}^r |w^i|^2.
		\end{equation}
		The proof of this fact follows after considering that $P_S=\sum_{i=1}^r s^i\otimes s^i=\sum_{i=1}^r w^i\otimes w^i$ and therefore their diagonal matrices coincide $\diag\left(P_S\right)
		=\diag\left(\sum_{i=1}^r s^i\otimes s^i\right)
		=\sum_{i=1}^r |s^i|^2 =\diag\left(\sum_{i=1}^r w^i\otimes w^i\right)= \sum_{i=1}^r |w^i|^2
		$.
\end{remark}	
The element
\begin{equation}
\label{def centroide}
c(m_{S,E})=\frac{1}{r} \sum_{i=1}^r |s^i|^2=\frac 1{\dim(S)}\diag(P_S)
\end{equation}
for (any) orthogonal basis $\{s^1,s^2,\dots,s^r\}$ of $S$ fulfils some interesting symmetric properties in the moment set $m_S$.
\begin{proposition}		\label{prop: varias props centroide}
			\begin{enumerate}
Let $S$ be a non trivial subspace of $\C^n$, $E$ a fixed basis of $\C^n$ and $c(m_S)=c(m_{S,E})$ defined as in \eqref{def centroide}. Then $c(m_S)$ satisfies the following properties.
	\begin{enumerate}
			\item $c(m_S)\in m_S$. 
			\item $c(m_S)$ coincides with the barycentre or centroid of the simplex generated by  $\{|w^1|^2,|w^2|^2,\dots,|w^r|^2\}\subset \R^n_{\geq 0}$ obtained from any orthonormal basis $\{w^1, w^2,\dots,w^r\}$ of $S$.
			\item Let $S$ and $V$ be subspaces of $\C^n$, with one of them not trivial, $S\perp V$, with dim$(S)=r$ and dim$(V)=k$. Then 
			\begin{equation}
			\label{eq ecuac centroide suma subesp perps}
			c\left(m_{S\operp V}\right)=\frac{1}{r+k}(r\, c(m_S)+k\, c(m_V)).
			\end{equation} This can be generalized to any number of mutually orthogonal subspaces. 
			\item $c(m_{\C^n})=(\frac{1}{n},\frac{1}{n},\dots,\frac{1}{n})$.
			\item If dim$(S)=r<n$, then $c(m_{S^\perp})=\frac{1}{n-r}( (1,1,\dots,1)- r\, c(m_S))$.
			\item If dim$(S)=r<n$ then the $i^\text{th}$ coordinate of $c(m_S)$ satisfies $c(m_S)_i\leq \frac1r$ for every $i=1,\dots, n$. 
			\item Given a subspace $D\subset S$, with dim$(D)=d<$ dim$(S)=r$, then $c(m_{S\ominus D})=c(m_{S\cap D^\perp})=\frac1{r-d}\left(r \, c(m_S)-d\, c(m_D)\right)$.
			\item  Let $S$ and $V$ be two subspaces of $\C^n$ with dimensions $r$ and $k$ respectively and $D= S \cap V$ of dimension $d$ such that $\left(S\cap D^{\perp}\right)\perp \left(V\cap D^{\perp}\right)$ holds. Then $c(m_{S+V})=\frac{1}{r+k-d}(r\, c(m_S)+k\, c(m_V) - d\, c(m_D))$.
					\end{enumerate}
	\end{enumerate}
\end{proposition}	
\begin{proof}
	\begin{enumerate}
		\item[(a) ] This follows after considering the definition of $c(m_S)$ \eqref{def centroide} and the statement \eqref{defequiv mS 3 en prop} in Proposition \ref{prop: equivalencias de momento}.  
		\item[(b) ] The centroid or barycentre of the simplex of $r$ vectors $\{w^i\}_{i=1}^r \subset\R^n$ is $\frac{1}{r} \sum_{i=1}^r w^i$. Then consider \eqref{igualdad_de_centroides_expandidos} and \eqref{def centroide}.
		\item[(c) ] This property can be proved using an orthonormal basis of $S\oplus V$ of elements of $S$ and $V$ and then \eqref{def centroide}.
		\item[(d) ] and (e) are apparent.
		\item[(f) ]  Using the previous (c), (d) and (e) follows that $c(m_{\C^n})_i=1/n=\frac{r\, c(m_S)_i+ (n-r) \, c(m_{S^\perp})_i}{n}$, which implies that $r\, c(m_S)_i=1-(n-r)\, c(m_{S^\perp})_i\leq 1$.
		\item[(g) ] Consider a completion of an orthonormal basis of $D$ to one of $S$ and use \eqref{def centroide}.
		\item[(h) ] This equation can be proved using similar techniques of basis completions and (f).
	\end{enumerate}
\end{proof}

\begin{remark}\label{rem comentario sobre centroide}
	Note the similarity of the equation \eqref{eq ecuac centroide suma subesp perps} with the one used to calculate the geometric centroid or barycentre of $m$ disjoint sets $A_j$ with $j=1,\dots,m$ using $c(\cup_{j=1}^m A_j)=\frac{\sum_{j=1}^m c(A_j) \mu(A_j)}{\sum_{j=1}^m\mu(A_j)}$ where $\mu$ is the corresponding measure.
\end{remark}
 
	\begin{proposition}
		Let $S$ and $W$ be two subspaces of $\C^n$ such that $S\perp W$ and such that there exist two corresponding orthonormal basis $\{s^h\}_{h=1}^r$ and $\{w^j\}_{j=1}^k$ that satisfy
		\begin{equation}\label{eq: condic S W soporte}
		\sum_{h=1}^{r} t_h|s^h|^2=\sum_{j=1}^{k} u_j |w^j|^2
		\end{equation}
		for $t_h, u_j\geq 0$ and $\sum_{h=1}^r t_h=\sum_{j=1}^k u_j=1$. That is, $(S,W)$ form a support (see Definition 1 and theorems 2 and 3 in \cite{soportes}).
		\\
		Then
		\begin{equation}\label{eq: coordenada de ligadura menor a 1/2}
		\left(\sum_{h=1}^{r} t_h|s^h|^2\right)_m = \left(\sum_{j=1}^{k} u_j |w^j|^2 \right)_m\leq \frac12,\ \ \text{ for every } m=1,\dots,n.
		\end{equation}
		\end{proposition}
		\begin{proof}
			Suppose there is $m_0$ such that 
			$
			\left(\sum_{h=1}^{r} t_h|s^h|^2\right)_{m_0} > \frac{1}{2}.
			$
			Then define the vectors $x=\sum_{h=1}^{r} \sqrt{t_h} \ e^{-i\, \text{arg}(s_{m_0}^h)} \ s^h$ and $y=\sum_{j=1}^{k} \sqrt{u_j} \ e^{-i\, \text{arg}(w_{m_0}^j)} \ w^j$.
			Observe that $x\in V$ and $\|x\|=\langle x, x\rangle^{1/2}  = \left(\sum _{h=1}^{r} t_h\right)^{1/2} =1$. 
			Similarly, follows that $y \in W$, $\|y\|=1$ and $x\perp y$. 
			\\
			Observe that the $m_0$ coordinate of $x$ satisfies  
			$
			|x_{m_0}|^2=\left(\sum_{h=1}^{r} \sqrt{t_h} \ e^{-i\, \text{arg}(s^h_{m_0})} \ s^h \right)_{m_0}^2 =
			\left(\sum_{h=1}^{r} \sqrt{t_h}\, |s^h_{m_0}| \right)^2 
			\geq\sum_{h=1}^{r} t_h |s^h_{m_0}|^2
			> \frac{1}{2}
			$, and in the case of $y$ also $|y_{m_0}|^2>\frac{1}{2}$.
			\\
			Now consider the dimension two subspace $C=\text{gen}\{x,y\}$ where $\{x,y\}$ is an orthonormal basis of $C$ and $c(m_C)=\frac{1}{2}(|x |^2 + |y |^2 )$ (see Proposition \ref{prop: varias props centroide} (c)). Then $c(m_C)_{m_0}= \frac{|x_{m_0}|^2 + |y_{m_0}|^2}{2}>  \frac{1/2 + 1/2}{2}=\frac{1}{2}$ which contradicts Proposition \ref{prop: varias props centroide} (f).
		\end{proof}
	\begin{remark}
		Note that the condition \eqref{eq: condic S W soporte} is equivalent to $m_S\cap m_W\neq \emptyset$ (see \eqref{defequiv mS 3 en prop} in Proposition \ref{prop: equivalencias de momento}). Then \eqref{eq: coordenada de ligadura menor a 1/2} implies that $m_S\cap m_W\subset [0,1/2]^n$ must always hold. This information is relevant since the supports $(S,W)$ allow the construction of minimal matrices (see Theorem 3, \cite{soportes}).
	\end{remark}


\section{Generic subspaces and its principal standard vectors}\label{Generic subspaces and its principal standard vectors}

\begin{definition}\label{def: subesp generico} 
	We call a subspace $S$ of $\C^n$ a \textbf{generic subspace} with respect to the basis $E=\{e_i\}_{i=1}^n$ if for every $j=1,\dots,n$ there is $x\in S$ such that $x_j=\langle x,e_i\rangle\neq 0$. 
	
	This condition is equivalent to any of the following statements: 
	\begin{itemize}
		\item $S$ is not included in a $(n-1)$-dimensional coordinate space $C_j=\{ z\in\C^n: z_j=0\}=(\text{span}\{e_j\})^\perp$ for $j=1,\dots,n$,
		\item or $\left(P_S(e_j)\right)_j\neq 0$ for all $j$ (otherwise $0=\langle P_S(e_j),e_j\rangle=\langle P_S(e_j),P_S(e_j)\rangle=\|P_S(e_j)\|^2$ which implies $e_j\perp S$).
	\end{itemize} 
\end{definition}

\begin{remark} 
	Observe that if $S$ is not a generic subspace of $\C^n$ we can find $m<n$ such that a natural immersion of $S$ in $\C^m$ becomes a generic subspace of $\C^m$.
\end{remark}

\begin{definition}\label{def: principal standard vectors} \textbf{Principal standard vectors.} 
		 If $P_S(e_j)=P_S e_j\neq 0$ we will denote with
	\begin{equation}\label{def: v a la j}
	v^j= \frac{P_S  e_j}{\|P_S   e_j\| }\in S , j=1, \dots, n
	\end{equation} 
	the unique principal vectors related to the standard basis $E=\{e_i\}_{i=1}^n$ that satisfy $(v^j)_j=v^j_j>0$ ($v^j_j=\langle\frac{P_S  e_j}{\|P_S   e_j\| },e_j\rangle=\|P_S e_j\|$)	
 and  minimize the angle between the subspace $S$ and the one dimensional coordinate axes span$(\{e_j\})$, that is
 \[
 \langle v^j,e_j\rangle=\max_{s\in S, \|s\|=1} |\langle s, e_j\rangle|.
 \]
 Note that the existence of these principal standard vectors for all $j$ requires that $S$ is a generic subspace (where $P_S e_j\neq 0$ for all $j$, see Definition \ref{def: subesp generico}). 
\end{definition}
To prove the uniqueness of $v^j$ suppose there exists $w\in S$ ($w\neq v^j$) with $\|w\|=1$ and such that 
\[
\langle w,e_j\rangle=\max_{s\in S, \|s\|=1} |\langle s, e_j\rangle|.
\]
Then  $w_j=v^j_j> 0$ and $v^j+w\neq 0$. If we suppose that $\|v^j+w\|=2$ then it can be proved that $|\langle v^j,w\rangle|\geq 1$ and therefore $v^j=\lambda w$ with $|\lambda|=1$, but since $v^j_j=w_j>0$ follows that $\lambda=1$ and then $v^j=w$. 
Now suppose $\|v^j+w\|<2$, and define $x=\frac{v^j+w}{\|v^j+w\|}\in S$, then
\begin{equation}
\label{eq: cuenta unicidad vj}
x_j=\langle x, e_j\rangle=\frac{\sfrac12(v^j_j+w_j)}{\sfrac12\|v^j+w\|}>1/2(v^j_j+w_j)=v^j_j,
\end{equation}
a contradiction.

\begin{remark}\label{rem: props vjk vkj etc}
We state here some properties of the principal standard vectors of a generic subspace $S$ defined in \eqref{def: v a la j} that follow after direct computations.
\begin{enumerate} \label{rem: props vec ppales}
	\item[a) ] For every $j=1, \dots, n$, $v^j_j=\|P_S(e_j)\|>0$.
	\item[b) ] If we also denote with $P_S$ the corresponding $n\times n$ matrix in standard coordinates, then its j$^\text{th}$-column satisfies $\text{col}_j(P_S)=P_Se_j= v^j_j v^j$ and its $j,k$ entry $(P_S)_{j,k} = v^j_j v^j_k$. In particular $(P_S)_{j,j}=\left(v^j_j\right)^2$.
	\item[c) ] Since the matrix $P_S$ is hermitian then $v^j_j v^j_k=\overline{v^k_k v^k_j}=v^k_k \overline{v^k_j}$ and then $\arg{(v^j_k)}=-\arg{(v^k_j)}$.
	\item[d) ] Directly from the previous item follows that if $P_Se_j$ and $P_Se_k$ are not null, then
	\begin{equation}
	\label{eq: igualdad fracc vjk/vjj etc}
	\frac{ v_k^j }{v_k^k} = \frac{ \overline{v_j^k }}{v_j^j} \ \text{ and }\ \frac{ \overline{v_k^j} }{v_k^k} = \frac{v_j^k }{v_j^j}.
	\end{equation}
	\item[e) ] 
				Observe that $ 0=v^j_k=\langle v^j,e_k\rangle=\langle v^j,P_S e_k\rangle=\|P_S(e_k)\| \langle v^j, v^k\rangle \Leftrightarrow v^j\perp v^k$ or $P_S(e_k)=0$  or $P_S(e_j)=0$. Therefore, if $S$ is  a generic subspace
				 \begin{equation}
				 \label{eq: vjk 0 sii vj perp vk}
				 v^j_k=0 \Leftrightarrow v^j\perp v^k.
				 \end{equation}

\end{enumerate}
\end{remark}
\begin{lemma}\label{lem: extremalidad de la coordenada}
	Let $S$ be a generic subspace of $\C^n$ and $w \in S$, with $\| w \|=1$. Then the following properties hold
	\begin{enumerate}
		\item 	$|v_j^j|=v_j^j\geq |w_j|$ for all $j=1, \cdots,n$.
		\item		Moreover,  
			\begin{equation}\label{eq: si coinciden las coord j con ppales los vec son multiplos}
			v_j^j = |w_j|\ \Leftrightarrow w= e^{i  \arg(w_j)} v^j\ ,
			\end{equation}	and in particular
		\item	\begin{equation}
				\label{eq: valajsubj igual valaksubj modulo equiv igual modulo}
				v^j_j=|v^k_j|\Leftrightarrow v^k=e^{i \arg(v^k_j)} v^j \Leftrightarrow v^j=e^{i \arg(v^j_k)} v^k \Leftrightarrow 
				 |v^j_i|=|v^k_i|,\forall i=1,\dots n.
				\end{equation}
		\item	As a consequence
				\begin{equation}
				\label{eq: vj vk LI sii vjj neq vkj}
				\{v^j,v^k\} \text{ is linearly independent } \Leftrightarrow v^j_j\neq |v^k_j| \ (v^j_j> |v^k_j|)
						\Leftrightarrow v^k_k\neq |v^j_k|\  (v^k_k> |v^j_k|).
				\end{equation}
		\end{enumerate}
\end{lemma}
\begin{proof} The first statement is apparent from the definition of $v^j$.
	For the second statement, observe that $v^j_j=\|P_S e_j\|> 0$ (see Remark \ref{rem: props vec ppales}). Then $w$, with $\|w\|=1$ is a multiple of $v^j$ if and only if $w=e^{i\arg(w_j)} v^j$. Then
\begin{equation}
\begin{split}
	w= {e}^{i \arg(w_j)}  v^j & \Leftrightarrow   |\langle w , v^j \rangle| = 1 \Leftrightarrow \left|\langle w , \frac{P_S  e_j}{\|P_S   e_j\| } \rangle\right| = 1 \Leftrightarrow \left|\langle P_S  w , \frac{e_j}{v_j^j} \rangle\right| = 1 \Leftrightarrow \\ 
& \Leftrightarrow  \frac{ |\langle w, e_j \rangle|}{v_j^j} = 1 \Leftrightarrow |w_j | = v_j^j
\end{split}
\end{equation}
The statement \eqref{eq: valajsubj igual valaksubj modulo equiv igual modulo} follows after replacing $w$ from \eqref{eq: si coinciden las coord j con ppales los vec son multiplos} with $v^k$ and applying some of the properties listed in Remark \ref{rem: props vjk vkj etc}.
\end{proof}

\begin{remark} \label{rem: coord j de vector en S}
	Given $S$ a generic subspace and $x \in S$ with  $\|x\|=1$, then its $j^\text{th}$-coordinate can be calculated as
	\begin{equation}\label{eq: vector por e_j}
	x_j= \| P_S   e_j \| \langle x \ , \ v^j \rangle=v^j_j  \langle x \ , \ v^j \rangle
	\end{equation} 
	This follows since
	$
	x_j=\langle x \ , \ e_j \rangle = \langle P_S   x \ , \ e_j \rangle = \langle x \ , \ P_S   e_j \rangle = 
	\langle x \ , \ \| P_S   e_j \|    v^j \rangle = \| P_S   e_j \| \langle x \ , \ v^j \rangle $. Therefore, $ x_j = 0$  if and only if $\langle x \ , \ v^j \rangle = 0$ or $\|P_S e_j \| =0 $, but this second condition cannot happen in a generic space.
	
\end{remark}

\begin{proposition}\label{prop: extremalidad de la proyeccion canonica}
	Given a subspace $S$ of $\C^n$ and $e_j$ a member of the standard basis $E$
	such that $P_S  e_j$ is not null (for example if $S$ is a generic subspace), then $|v^j|^2=(|v^j_1|^2, |v^j_2|^2,\dots, |v^j_n|^2) $ 
	is an extreme point in $m_S$. 
	
	Moreover, if $|v^j|^2 $ is a convex combination of $|y|^2$ and $|z|^2$ with $y, z\in S$, then $y$ and $z$ must be multiples of $v^j$.
	
	\begin{proof}  
		Suppose that there are two vectors $y, z \in S$, with $\|y\|= \|z \| =1$ such that $|v^j|^2 = t |y|^2 + (1-t)|z|^2$ with $0 \leq t \leq 1$.
		Then in particular $|v^j_j|^2 = t |y_j|^2 + (1-t)|z_j|^2$ which implies that  $|v^j_j|=v^j_j=|y^j_j|=|z^j_j|$. Then using \eqref{eq: si coinciden las coord j con ppales los vec son multiplos} in Lemma \ref{lem: extremalidad de la coordenada} we obtain that both $y$ and $z$ must be multiples (by a complex with modulus one) of $v^j$. This implies that $|v^j|^2=|y|^2=|z|^2$.
	\end{proof}	
\end{proposition}

\section{Curves of extreme points in $m_S$}\label{curves of extreme points}
\begin{definition}\label{def: jflechakV}
Let $S$ be a generic subspace and $v^j$ and $v^k$ (defined in \ref{def: v a la j}) be linearly independent. We define a curve $\dejak{j}{v}{k}:[0,\pi/2]\to S$ that starts in $v^j$ and passes through $e^{i\arg(v^j_k)}v^k$:
\begin{equation}
\label{eq: def vtildejk}
\dejak{j}{v}{k}(t)=\cos(t)v^j+\sin(t)\ e^{i \arg\left(v^j_k\right)} \frac{\left(v^k -\langle v^k,v^j\rangle v^j\right)}{\|v^k -\langle v^k,v^j\rangle v^j\|}  , \ t\in[0,\pi/2].
\end{equation}
\end{definition}
Note that $\big\|\dejak{j}{v}{k}(t)\big\|=1$ and that using basic properties of $v^j$ and $v^k$ (see Remark \ref{rem: props vec ppales}) the $j$ and $k$ coordinates of $\dejak{j}{v}{k}(t)$ (that we denote with $\dejak{j}{v_j}{k}(t)$ and $\dejak{j}{v_k}{k}(t)$) can be computed
\begin{equation}
\label{eq: props curva vjk}
{ \dejak{j}{v_j}{k}}(t)=\cos(t) v^j_j 
\ ,\ \text{  and  } \  \ 
\dejak{j}{v_k}{k}(t) = \cos(t) v^j_k+\sin(t) e^{i \arg\left(v^j_k\right)} \sqrt{(v^k_k)^2-|v^j_k|^2}.
\end{equation}

\begin{remark}
	The curve vectors $\dejak{j}{v}{k}(t)$ can be obtained by projections of particular linear combinations of the $e_j$ and $e_k$ standard vectors. Namely,
	\begin{equation}\label{eq: vjk como proy de curva en ej ek}
	\dejak{j}{v}{k}(t)=P_S\left(\dejak{j}{\beta}{k}(t)\right), \ \text{ for } \dejak{j}{\beta}{k}(t)=	
	 \cos(t)\frac{e_j}{\|P_S(e_j)\|}+\sin(t)\ e^{i \arg\left(v^j_k\right)} \frac{\left(\frac{e_k}{\|P_S(e_k)\|} -\langle v^k,v^j\rangle \frac{e_j}{\|P_S(e_j)\|}\right)}{\|v^k -\langle v^k,v^j\rangle v^j\|}.
	\end{equation}
	This implies that for each $t\in[0,\pi/2]$, the vector $\dejak{j}{v}{k}(t)$ attains the minimal angle between $S$ and the one dimensional subspace spanned by $\dejak{j}{e}{k}(t)=\frac{\dejak{j}{\beta}{k}(t)}{\big\|\dejak{j}{\beta}{k}(t)\big\|}$   
	\[
	\langle\dejak{j}{v}{k}(t),\dejak{j}{e}{k}(t)\rangle=\max_{s\in S, \|s\|=1}	|\langle s,\dejak{j}{e}{k}(t)\rangle|.
	\]
	Note that $\dejak{j}{e}{k}(t)$ is included in span$\{e_j, e^{i \arg\left(v^j_k\right)}\ e_k\}= $ span$\{e_j, e_k\}$ for all $t\in[0,\pi/2]$.
	
	Moreover, it can be computed that $\langle\dejak{j}{v}{k}(t),\dejak{j}{e}{k}(t)\rangle=\left\|P_S\left(\dejak{j}{e}{k}(t)\right)\right\|=
	\frac{1}{\big\|\dejak{j}{\beta}{k}(t)\big\|}>0$ (if zero then $v^j$ and $v^k$ should be linearly dependent), $\dejak{j}{v}{k}(t)=\frac {P_S\left(\dejak{j}{e}{k}(t)\right)}{\big\|P_S\left(\dejak{j}{e}{k}(t)\right)\big\|}$, 
	and $\dejak{j}{v}{k}(t)$ is unique among the unit vectors $s\in S$ that attain this maximum with the property that $\langle s,\dejak{j}{e}{k}(t)\rangle >0$ (this can be proved as done in \eqref{eq: cuenta unicidad vj}).
\end{remark}

Following the same procedures as those used in Lemma  \ref{lem: extremalidad de la coordenada} for $v^j$, similar results can be obtained for $\dejak{j}{v}{k}(t)$ as stated in the next lemma.
		\begin{lemma}\label{lem: extremalidad de la coordenada ejk}
			Let $S$ be a generic subspace of $\C^n$ and $w \in S$, with $\| w \|=1$. Then the following properties hold
			\begin{enumerate}
				\item 	$|\langle \dejak{j}{v}{k}(t), \dejak{j}{e}{k}(t)\rangle|=\langle \dejak{j}{v}{k}(t), \dejak{j}{e}{k}(t)\rangle \geq |\langle w, \dejak{j}{e}{k}(t)\rangle|$ for all $t\in[0,\pi/2]$.
				\item		Moreover,  
				\begin{equation}\label{eq: si coinciden las coord jk con ppales los vec son multiplos}
				\langle \dejak{j}{v}{k}(t), \dejak{j}{e}{k}(t)\rangle = |\langle w, \dejak{j}{e}{k}(t)\rangle| 
				\ \Leftrightarrow 
				w= e^{i  \arg(\langle w, \dejak{j}{e}{k}(t))\rangle}\  \dejak{j}{v}{k}(t)\ ,
				\end{equation}	
				\item		
				and in particular
				\begin{equation}
				\label{eq: vjkt igual vjks modulo equiv igual modulo}
				\begin{split}
				\langle \dejak{j}{v}{k}(t),& \dejak{j}{e}{k}(t)\rangle=
				|\langle \dejak{j}{v}{k}(s_0), \dejak{j}{e}{k}(t)\rangle|\ \text{, for } s_0\in[0,\pi/2]
				\Leftrightarrow \\
				&\Leftrightarrow 
				\dejak{j}{v}{k}(s_0)= e^{i  \arg(\langle \dejak{j}{v}{k}(s_0), \dejak{j}{e}{k}(t))\rangle}\  \dejak{j}{v}{k}(t)  \\
				&\Leftrightarrow  
				|\langle \dejak{j}{v}{k}(t), \dejak{j}{e}{k}(u)\rangle|=
				|\langle \dejak{j}{v}{k}(s_0), \dejak{j}{e}{k}(u)\rangle|\ ,\forall u\in[0,\pi/2].
				\end{split}
				\end{equation}
				\item	As a consequence
				\begin{equation}
				\label{eq: vjkt vjks LI sii vjkt neq vjks}
				\left\{\dejak{j}{v}{k}(t) , \dejak{j}{v}{k}(s)\right\} \text{ is linearly independent } \Leftrightarrow \langle \dejak{j}{v}{k}(t), \dejak{j}{e}{k}(t)\rangle\neq
				|\langle \dejak{j}{v}{k}(s), \dejak{j}{e}{k}(t)\rangle|.
				\end{equation}
			\end{enumerate}
		\end{lemma}
	
\begin{proposition}\label{propo: props curva vjk} Let $S$ be a generic subspace of $\C^n$ and $v^j$ and $v^k$ from \eqref{def: v a la j}. Consider $\overset{j \curvearrowright k}{v} : [0, \frac{\pi}{2}] \rightarrow \text{Im}\big(\overset{j \curvearrowright k}{v}\big) \subset S$  the curve defined in \eqref{eq: def vtildejk}, then:
\begin{enumerate}
	\item the map $\overset{j \curvearrowright k}{v}$ is bijective,
	\item for all $t \in [0, \frac{\pi}{2}]$, holds that
	 \begin{equation}\label{producto escalar entre v^t y v^j, v^k}
		\langle \overset{j \curvearrowright k}{v}(t) , v^j \rangle \in  \R_{\geq 0}  \qquad 
		\text{and} \quad \langle \overset{j \curvearrowright k}{v}(t) , \textit{e}^{\textit{i} \arg\left(v_k^j\right)}  v^k  \rangle \in  \R_{\geq 0}
		\end{equation} 
   \item $\overset{j \curvearrowright k}{v}(0) = v^j$ \quad and \quad $\dejak{j}{v}{k}\left( \arccos\left( { |v_k^j| }/{v_k^k}\right) \right)= e^{\textit{i} \arg(v_k^j)}  v^k $.
 \end{enumerate}
\end{proposition}
	\begin{proof}	
		\begin{enumerate}
			\item To prove the bijectivity of $\overset{j \curvearrowright k}{v}$ it is enough to observe that the $j$-coordinate of the curve is $\overset{j \curvearrowright k}{v_j}(t) = \cos(t) v_j^j $ (see \eqref{eq: props curva vjk}) which is a strictly decreasing real function from $[0,\frac{\pi}{2}]$ onto $[0, v_j^j]$.
			\item Fix $t\in[0,\frac{\pi}{2}]$ and recall that $\overset{j \curvearrowright k}{v}(t) $ has norm 1 (see \eqref{eq: def vtildejk}). Moreover, the fact that $S$ is a generic subspace implies that $v_j^j, v_k^k \in \R_{>0}$,
			$e^{i \arg(v_k^j)}  v^k_j\in\R_\geq 0$ and $\arg(v_k^j)
			=\arg(e^{\textit{i} \arg(v_k^j)}  v^k_k)$.
			
			On the other hand 
			\[
			\arg(\langle v^j, \textit{e}^{\textit{i} \arg(v_k^j)}  v^k \rangle) = 
			\arg(\textit{e}^{-\textit{i} \arg(v_k^j)}  \langle v^j,  v^k \rangle) = \arg\left(\|P_S e_k \| \langle v^j,  e_k \rangle\right) -\arg(v_k^j)= 0	
			\]
			and consequently
			\begin{equation}
			\label{argumentos de v^t_(1,2)}
			v_j^j,\, e^{i \arg(v_k^j)}  v^k_j,\,  e^{i \arg\left(v^j_k\right)} \frac{\left(v^k_j -\langle v^k,v^j\rangle v^j_j\right)}{\|v^k -\langle v^k,v^j\rangle v^j\|} \in \R_{\geq 0} \text{ , hence }
			\ \overset{j \curvearrowright k}{v_j}(t) \in \R_{\geq 0} .
			\end{equation}
			Similarly
			\begin{equation}
			\label{eq: igualdad arg vjk con arg jvkkt}
						\qquad \arg(v_k^j) = \arg( e^{i \arg(v_k^j)}  v^k_k ) = \arg\left(e^{i \arg\left(v^j_k\right)} \frac{\left(v^k_k -\langle v^k,v^j\rangle v^j_k\right)}{\|v^k -\langle v^k,v^j\rangle v^j\|}\right) = \arg\left(\overset{j \curvearrowright k}{v_k}(t)\right)
			\end{equation}
				The previous equalities \eqref{argumentos de v^t_(1,2)} and \eqref{eq: igualdad arg vjk con arg jvkkt} and Remark \ref{rem: coord j de vector en S} imply that
			\[
			\langle \overset{j \curvearrowright k}{v}(t) , v^j \rangle  = \overset{j \curvearrowright k}{v_j}(t) {/v^j_j} \in  \R_{\geq 0}  
			\]
			and
			\[
			 \langle \overset{j \curvearrowright k}{v}(t) , \textit{e}^{\textit{i}  \arg(v_k^j)}  v^k  \rangle = \textit{e}^{ - \textit{i}   \arg(v_k^j)} \langle \overset{j \curvearrowright k}{v}(t) ,  v^k  \rangle 
			\overset{\eqref{eq: igualdad arg vjk con arg jvkkt}}{=} e^{ - i \arg\left( \overset{j \curvearrowright k}{v_k}(t) \right)}\   \overset{j \curvearrowright k}{v_k}(t) {/v^k_k} \in \R_{\geq 0}. 
			\]
			
			\item 
			
			It is trivial that $\overset{j \curvearrowright k}{v}(0) = v^j$.
			Now if $v^j$, $v^k$ are linearly independent, then $0<\|v^k-\langle v_k,v^j\rangle v^j\|=\sqrt{1-\frac{|v^k_j|^2}{(v_j^j)^2}}=\frac{\sqrt{(v_j^j)^2-|v^k_j|^2}}{v^j_j}<1$, and therefore there exists $t_0\in (0,\pi/2)$ such that
			\[
			\sin(t_0)=\|v^k-\langle v_k,v^j\rangle v^j\|=\frac{\sqrt{(v_j^j)^2-|v^k_j|^2}}{v^j_j}. 
			\]
			Using that
			$
			\left(\frac{\sqrt{(v_j^j)^2-|v^k_j|^2}}{v^j_j}\right)^2+ \left(|v^k_j|/v^j_j\right)^2=1
			$, then 
			$
			t_0=\arccos(|v^k_j|/v^j_j)=\arccos(|v^j_k|/v^k_k)\in (0,\pi/2)
			$ 
		since $|v^k_j|/v^j_j=|v^j_k|/v^k_k<1$. Otherwise, if $|v^k_j|/v^j_j=1$, $v^j$ and $v^k$ must be linearly dependent (see \eqref{eq: vj vk LI sii vjj neq vkj} in Lemma \ref{lem: extremalidad de la coordenada}).
			 \\
			Evaluating the sine and cosine in that $t_0$ we obtain:
			\begin{equation}
			\begin{split}
			\dejak{j}{v}{k}\left(t_0\right)& =  \cos(t_0) v^j+\sin(t_0) \ e^{i \arg\left(v^j_k\right)} \frac{\left(v^k -\langle v^k,v^j\rangle v^j\right)}{\|v^k -\langle v^k,v^j\rangle v^j\|}\\
			& = \left(|v_k^j|/{v_k^k}\right) v^j+ { \|v^k -\langle v^k,v^j\rangle v^j\|} \ e^{i \arg\left(v^j_k\right)} \frac{\left(v^k -\langle v^k,v^j\rangle v^j\right)}{ {\|v^k -\langle v^k,v^j\rangle v^j\|}}\\
			&= \left(|v_k^j|/{v_k^k}\right) v^j+ \, e^{i \arg\left(v^j_k\right)}  {\left(v^k -\langle v^k,v^j\rangle v^j\right)} 
			\end{split}
			\end{equation}
			And using that $e^{i\arg(v^j_k)}\langle v^k,v^j\rangle = e^{i\arg(v^j_k)} v^k_j/v^j_j= {|v^k_j|}/{v^j_j}$ we obtain
			\begin{equation}
			\begin{split}
			\dejak{j}{v}{k}\left(t_0\right)&=\left(|v_k^j|/{v_k^k}\right) v^j+  e^{i \arg\left(v^j_k\right)}  \left(v^k -\langle v^k,v^j\rangle v^j\right)= {\left(|v_k^j|/{v_k^k}\right) v^j}+  e^{i \arg\left(v^j_k\right)}  v^k - ({|v^k_j|}/{v^j_j}) v^j \\
			&=  e^{i \arg\left(v^j_k\right)}  v^k.
			\end{split}
			\end{equation}

		\end{enumerate}
	\end{proof}

\begin{proposition}\label{prop: desigs entre coords vj xj vk xk} Consider the vectors $\dejak{j}{v}{k}(t)$ from Definition \ref{def: jflechakV}, $x\in S$ with $\|x\|=1$ and $\overline{e_je_k}$ the segment between $e_j$ and $e_k$ projected to the $j$ and $k$ coordinates (in $\R^2$). 
	Then there exists $t_x\in [0,\pi/2]$ such that
\begin{equation}
\label{eq: desigualdades coord j y k}
	 |x_j|\leq\big|\dejak{j}{v_j}{k}(t_x) \big|\ \ ,\ \ |x_k|\leq \big|\dejak{j}{v_k}{k}(t_x)\big|\ \text{ and}
\end{equation}	
\begin{equation}
\label{eq: curva posta esta mas cerca de segmento}
\operatorname{dist}\left(\big(\big|\dejak{j}{v_j}{k}(t_x)\big|^2,\big|\dejak{j}{v_k}{k}(t_x)\big|^2\big),\overline{e_j e_k}\right) \leq \dist\left((|x_j|^2,|x_k|^2),\overline{e_j e_k}\right).
\end{equation} 
\end{proposition}
\begin{proof}	
			 In order to alleviate the notation let us write 
			 \begin{equation}
			 \label{notacion wjkMonio}
			 \widetilde{w}^{jk}=e^{i \arg\left(v^j_k\right)} \frac{\left(v^k -\langle v^k,v^j\rangle v^j\right)}{\|v^k -\langle v^k,v^j\rangle v^j\|}.
			 \end{equation} 
			 Recall that $ \dejak{j}{v_j}{k}(t)=\cos(t) v^j+\sin(t) \wjk$ with $v^j\perp \wjk$ (see \eqref{eq: def vtildejk}).
			 Then we can write $x\in S$ as a linear combination of $v^j,\, \widetilde{w}^{jk}$ and a vector $y\in S$ with $\|y\|=1$, orthogonal to the subspace spanned by the other two vectors 
	 \begin{equation}
	 \label{eq: descomp x = a b c}
	 	 x = a\ v^j + b\ \widetilde{w}^{jk} + c\ y \text{ , with } a, b , c \in \C.
	 \end{equation}
	 We will consider two cases: $c=0$ and $c\neq 0$.
	\begin{itemize}
		\item Case $c=0$
		
Recalling the $j^\text{th}$ and $k^\text{th}$ standard coordinates of $v^j$ and $\widetilde{w}^{jk}$  as in 
\eqref{eq: props curva vjk} we can write (suppose $j<k$):
 \[
 x = a\, v^j + b\, \widetilde{w}^{jk} =\left(\dots\,, a v^j_j\, ,\dots\, ,\, a v^j_k+b e^{\text{i} \arg(v^j_k)} \sqrt{(v^k_k)^2-|v^j_k|^2},\, \dots\right). 
\]
 The orthonormality of $\{v^j, \widetilde{w}^{jk}\}$ and $\|x\|=1$ implies that $|a|^2+|b|^2=1$.
Now define $\alpha=\arccos(|a|)\in[0,\pi/2]$ and consider 
\[
\dejak{j}{v}{k}(\alpha)=\cos(\alpha)\, v^j+ \sin(\alpha)\, \widetilde{w}^{jk}.
\]
Then
\begin{equation}
\label{eq: igualdad mod xj vjj}
|x_j|=|a| \, |v^j_j|=\big|\dejak{j}{v_j}{k}(\alpha) \big|
\end{equation}
and
\begin{equation}
\label{eq: xk menor o igual a curva en t_x Caso 1}
\begin{split}
|x_k|&=\left|a\, v^j_k +b\, e^{i \arg(v^j_k)} \sqrt{(v^k_k)^2-|v^j_k|^2}\right|\\
&\leq 
|a|\, |v^j_k| +|b|\, \left|\sqrt{(v^k_k)^2-|v^j_k|^2}\right|
= \cos(\alpha)  |v^j_k| + \sin(\alpha) \sqrt{(v^k_k)^2-|v^j_k|^2}\\
&= \big|\dejak{j}{v_k}{k}(\alpha)\big|
\end{split}
\end{equation}
where in the last equality we have used that $v^j$ and  $\widetilde{w}^{jk}$ have the same argument in its $k^\text{th}$ coordinate (see \eqref{eq: props curva vjk}).
Therefore, $|x_k|\leq \big|\dejak{j}{v_k}{k}(\alpha)\big|$ which together with \eqref{eq: igualdad mod xj vjj} proves \eqref{eq: desigualdades coord j y k} in this case. Then follows that
\[
\dist\left((\big|\dejak{j}{v_j}{k}(t_{x}) \big|^2,\big|\dejak{j}{v_k}{k}(t_{x})\big|^2),\overline{e_j e_k}\right) 
\leq 
\dist\left((|x_1|^2,|x_2|^2),\overline{e_j e_k}\right).  
\]

	\item Case $c \neq 0$ (and $y\neq 0$)
		
		Using Remark \ref{rem: coord j de vector en S} and the fact that $y \perp \text{span}\{v^j , \widetilde{w}^{jk} \}=\text{span}\{v^j , v^k \}$ it can be proved that $y_j = \|P_S e_j\| \langle y \ , \ v^j\rangle =0 $ and similarly  $y_k=0$. Then,		
		$
		x_j=a\, v_j^j + b\, \widetilde{w}^{jk}_j$  and  $x_k=a\,  v_k^j + b\,  \widetilde{w}^{jk}_k
		$.
		
		If $x_j=0$ \ and \ $x_k=0$ it is enough to take $t_x=0$ ($\dejak{j}{v}{k}(t_x)=v^j$). 
		
		Consider now the case when $x_j\neq 0$ or $x_k \neq 0$ and the define the vector $\widehat{x } \in S$ as
		\[
		\widehat{x} = \frac{P_{\text{span}\{v^j,\widetilde{w}^{jk}\}}(x)}{\|P_{\text{span}\{v^j,\widetilde{w}^{jk}\}}(x)\|}= \frac{a v^j+b \widetilde{w}^{jk}}{\| a v^j+b \widetilde{w}^{jk} \|}.
		\]
		Since $\| a v^j+b \widetilde{w}^{jk} \| < \|x\|=1$ 
		 (because $y\perp v^j$, $y\perp \widetilde{w}^{jk}$) and the $j^\text{th}$ and $k^\text{th}$ coordinates of $a v^j+b \widetilde{w}^{jk}$ and $x$ coincide ($y_j=y_k=0$), we have that 
		\begin{equation}
		\label{eq: desig coord j,k de xhat y x}
		|\widehat{x}_j|   \geq |x_j| \ \text{ and }\ |\widehat{x}_k|  \geq  |x_k|.
		\end{equation}
		 Since $\widehat{x}$ is included in the case already considered when $c=0$, there exists $t_{\widehat{x}}$ that satisfies \eqref{eq: desigualdades coord j y k}, (see \eqref{eq: igualdad mod xj vjj} and \eqref{eq: xk menor o igual a curva en t_x Caso 1}). That is, $|\widehat{x}_j|=\big|\dejak{j}{v_j}{k}(t_{\widehat{x}}) \big|$ and $|\widehat{x}_k|\leq \big|\dejak{j}{v_k}{k}(t_{\widehat{x}})\big|$. Then using \eqref{eq: desig coord j,k de xhat y x} we obtain
		\[
		|x_j|\leq\big|\dejak{j}{v_j}{k}(t_{\widehat{x}}) \big|\ \ ,\ \ |x_k|\leq \big|\dejak{j}{v_k}{k}(t_{\widehat{x}})\big|
		\]
		which in turn proves
		\[
		d\left((\big|\dejak{j}{v_j}{k}(t_{\widehat{x}}) \big|^2,\big|\dejak{j}{v_k}{k}(t_{\widehat{x}})\big|^2),\overline{e_j e_k}\right) \leq d\left((|x_1|^2,|x_2|^2),\overline{e_j e_k}\right).  
		\]
	\end{itemize}		
	\end{proof}

\begin{remark} There are many different $x\in S$ with $\|x\|=1$ that satisfy \eqref{eq: desigualdades coord j y k} for the same value of $t_x$ (see for example Figure \ref{fig: curva proyectada solo con modulos}). 
\end{remark}
%
 
  \begin{theorem}\label{teo: extremalidad de los puntos de la curva}
Consider the curve $\dejak{j}{v}{k}:[0,\pi/2]\to \R^n_{\geq 0}$ from Definition \ref{def: jflechakV} and $x\in S$, with $\|x\|=1$. 
Then there exists a unique $t_x \in [0,\frac{\pi }{2}]$ such that
 	\begin{equation}
 		\label{eq: igualdad y desig de xj xk en mom} 
 		|x_j|=|\dejak{j}{v_j}{k}(t_x)|  \quad \text{ and } 
 		\quad  |x_k|\leq |\dejak{j}{v_k}{k}(t_x)|.
 	\end{equation}
 	Moreover, if $x=a\, v^j+b\, \wjk +c\, y$ (with $\wjk$ as in  \eqref{notacion wjkMonio} and $y\perp v^j,  \wjk$), then $t_x=\arccos(|a|)$.
 \end{theorem}
  \begin{proof} Consider first the the existence of $t_x$.
 
  	We continue using the notation $\wjk$ as in  \eqref{notacion wjkMonio}.
 	As in \ref{eq: descomp x = a b c}, we write
 	\[
 	x = a\, v^j + b\, \widetilde{w}^{jk} + c\, y \text{ , with } a, b , c \in \C \text{ and } y\perp v^j, \wjk.
 	\] 
 	Now we consider different cases. If $c=0$ we can choose $ t_x = \arccos(|a|)$ as in \eqref{eq: igualdad mod xj vjj} and \eqref{eq: xk menor o igual a curva en t_x Caso 1} in the proof of the previous Proposition \ref{prop: desigs entre coords vj xj vk xk} to obtain the equality and inequality required. 
 	
 	Now consider the case $c\neq 0$ and some sub-cases
 	\begin{itemize}
 		\item If $a=0$ then it must be $x_j= 0 \ v^j_j + b\ \wjk_j+c\ 0=0$ where we have used that $\wjk_j=0$ (a direct computation) and $y_j=0$ (see Case $c\neq 0$ in Proposition \ref{prop: desigs entre coords vj xj vk xk}). Then we can choose
 		$t_x = \frac{\pi}{2}$, and obtain 	$|\dejak{j}{v_j}{k}(\pi/2)| =0 = |x_j|$.
 		
 		For the $k^\text{th}$-coordinate, $y_k=0$ and \eqref{eq: props curva vjk} imply
 		
 		$|x_k|=
 		\left|0\, v^j_k +b\, e^{i \arg(v^j_k)} \sqrt{(v^k_k)^2-|v^j_k|^2}\right|=
 		\left| b\, \sqrt{(v^k_k)^2-|v^j_k|^2} \right| \leq 
 		| 1 \, \wjk_k | =  |\dejak{j}{v_k}{k}(\frac{\pi }{2})|$  
 		\item Consider now $a \neq 0$. 
 		If we choose $t_x = \arccos(|a|)$, then	using again that $\wjk_j=0$ and $y_j=0$ we obtain
 		\[
 		|x_j| = | a \, v_j^j | = \cos(t_x) v_j^j= |\dejak{j}{v_j}{k}(t_x)|.
 		\] 		
 		Moreover, since $|a|^2 + |b|^2 + |c|^2=1$,
 		\begin{equation*}
 		\begin{split}
 		|x_k|&=
 		\left|a\, v^j_k +b\, e^{i \arg(v^j_k)} \sqrt{(v^k_k)^2-|v^j_k|^2}\right| \\
 		&\leq 
 		\left|a\, v^j_k + \sqrt{1 -|a|^2}\, e^{i \arg(v^j_k)} \  \sqrt{(v^k_k)^2-|v^j_k|^2}\right|= |\dejak{j}{v_k}{k}(t_x)|
 		\end{split}
 		\end{equation*}
 	\end{itemize}
 	which ends the proof of the existence of $t_x$.
 	\\ 
 	The uniqueness of $t_x\in[0,\pi/2]$ can be proved using (as in the proof of (1) in Proposition \ref{propo: props curva vjk}) the bijectivity of the map $t\mapsto \overset{j \curvearrowright k}{v_j}(t) = \cos(t) v_j^j$  for $t \in[0,\pi/2]$.
 	
\end{proof}
\begin{remark}
	The maximality of $v^i$ in its $i$-coordinate $v^i_i$ (see Lemma \ref{lem: extremalidad de la coordenada}) implies that if $x\in S$ with $\|x\|=1$ then $|x_j|\in[0,v^j_j]$ and $|x_k|\in[0,v^k_k]$. Then, the projection of  $m_S$ to the $j$ and $k$ coordinates is included in the rectangle $[0,(v^j_j)^2]\times [0,(v^k_k)^2]$. In the following results we will show more precise boundaries of $m_S$.
\end{remark}
\begin{theorem}\label{teo props proy vjk en coords jk}
	Let $S\subset\C^n$ be a generic subspace, $\{v^j,v^k\}$ two linearly independent principal standard vectors, $m_S$ the moment of $S$ as in Definitions \ref{def: principal standard vectors} and \ref{def: momento de un subespacio} respectively, and $\gamma_{j,k}:[0,\pi/2]:\to m_S\subset \R^n_{\geq 0}$ the curve defined by 
	\begin{equation}
	\label{def gamajk} 
	\gamma_{j,k}(t)=\big|\dejak{j}{v}{k}(t)\big|^2=\left(\big|\dejak{j}{v_1}{k}(t)\big|^2,\dots,\big|\dejak{j}{v_n}{k}(t)\big|^2\right)
	\end{equation}
	with $\dejak{j}{v}{k}(t)$ as in Definition \ref{def: jflechakV}. 
	
	Then 
	
	\begin{enumerate}
		\item the projection of the $j^\text{th}$ and $k^\text{th}$ coordinates of $\sqrt{\gamma_{j,k}}$ to $\R^2$ given by $t\mapsto \left(|\dejak{j}{v_j}{k}(t)|,|\dejak{j}{v_k}{k}(t)|\right)$ is part of an ellipse centred at the origin,
		
		\item if $v^j$ and $v^k$ are not orthogonal, the points  $\left(|\dejak{j}{v_j}{k}(t)|^2,|\dejak{j}{v_k}{k}(t)|^2\right)$ from
		the projected curve $\gamma_{j,k}$ to the plane spanned by $e_j$, $e_k$ are all extreme points of the projection of $m_S$ to the same plane,
		
		\item in case $v^j\perp v^k$ then $\left(|\dejak{j}{v_j}{k}(t)|^2,|\dejak{j}{v_k}{k}(t)|^2\right)$ parametrizes a segment that is in the boundary of the projection of $m_S$ to the plane spanned by $e_j, e_k$.
	\end{enumerate}
	%
	%
\end{theorem}

\begin{proof}	
	Let $P_{j,k}:\R^n\to\R^2$ be the projection $P_{j,k}(x)=(x_j,x_k)$, then the curve obtained in the plane after the projection of $\sqrt{\gamma_{j,k}(t)}$ is (see \eqref{eq: props curva vjk}) 
	\begin{equation}\label{eq: curva proyec mod j y k}
	\begin{split}
	P_{j,k}\left(\big|\dejak{j}{v}{k}(t)\big| \right)&=
	\left(\big|\cos(t) v^j_j \big|  \ ,\ \left| \cos(t) v^j_k+\sin(t) e^{i \arg\left(v^j_k\right)} \sqrt{(v^k_k)^2-|v^j_k|^2}\right| \right)\\
	&= \left(\cos(t) v^j_j   \ ,\ \cos(t) |v^j_k|+\sin(t) \sqrt{(v^k_k)^2-|v^j_k|^2}\right)\\
	&= \cos(t)\, \left(v^j_j, |v^j_k|\right) + \sin(t)\, \left(0,\sqrt{(v^k_k)^2-|v^j_k|^2} \right) 
	\end{split}
	\end{equation}	
	which is clearly part of an ellipse centred at the origin (recall that $t\in[0,\pi/2]$). It starts in $P_{j,k}(|v^j|)$, passes through $P_{j,k}(|v^k|)$ (see Proposition \ref{propo: props curva vjk} part (3)) and ends in the span$\{(0,1)\}$ axis. See for example Figure \ref{fig: curva proyectada solo con modulos}.
	
	\begin{figure}
		\begin{center}
			\epsfig{file=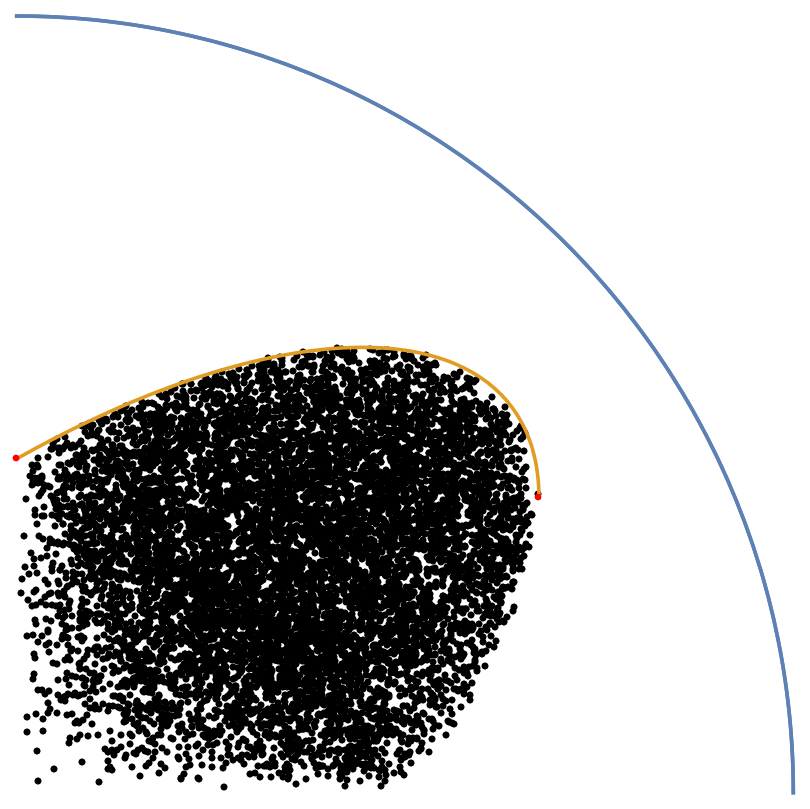, width=5cm}
		\end{center} 
		\caption{For $j=1$ and $k=2$ image of the curve $|\dejak{j}{v}{k}|$ projected on the $j$ and $k$ coordinates for a 3 dimensional subspace $S$ of $\C^9$ (in orange). The blue curve is a quarter of the unit circle and the black dots are projected points of $m_S$ (with the square roots of its entries) taken randomly.}
		\label{fig: curva proyectada solo con modulos}
	\end{figure}
	
	Now squaring the coordinates of \eqref{eq: curva proyec mod j y k} 	we obtain
	\begin{equation} \label{eq: segmento mas vector perp}
	\begin{split}
	P_{j,k}\left(\big|\dejak{j}{v}{k}(t)\big| \right)^2&=\bigg(\cos^2(t) (v^j_j)^2   \ ,\ \cos^2(t) |v^j_k|^2+\sin^2(t) \left((v^k_k)^2-|v^j_k|^2\right)+\\
	&\hskip3.5cm +
	2 \cos(t) |v^j_k| \sin(t) \sqrt{(v^k_k)^2-|v^j_k|^2}\bigg)\\
	&=\cos^2(t) \left( (v^j_j)^2   \ ,\  |v^j_k|^2\right)+ \sin^2(t)\left(0, \left((v^k_k)^2-|v^j_k|^2\right)\right)+\\
	&\hskip3.5cm +\left(0,
	2 \cos(t) |v^j_k| \sin(t) \sqrt{(v^k_k)^2-|v^j_k|^2}\right).
	\end{split}
	\end{equation} 
 This is the parametrization of a segment that joins $\left( (v^j_j)^2   \ ,\  |v^j_k|^2\right)$ with $\left(0, \left((v^k_k)^2-|v^j_k|^2\right)\right)$ plus a vertical vector with second coordinate $\geq 0$. Note that this second coordinate is zero only if $t=0$, $t=\pi/2$, $v^j_k=0$ or $v^k_k=|v^j_k|$. This last case is not possible because otherwise using \eqref{eq: vj vk LI sii vjj neq vkj} the vectors $v^j$ and $v^k$ would be linearly dependent and we are supposing in the hypothesis that they are not (see Definition \ref{def: jflechakV}). Then the curve in \eqref{eq: segmento mas vector perp} is a segment only when $v^j_k=0$, which is equivalent to
	\begin{equation}
	\label{eq: vjk cero sii vj perp vk}0=v^j_k=\langle v^j,e_k\rangle=\|P_S(e_k)\| \langle v^j, v^k\rangle \Leftrightarrow v^j\perp v^k
	\end{equation}
	\begin{enumerate}
		\item Case $v^j{\notperp} v^k$ (equivalently $v^j_k\neq 0$)
		
		The curve parametrized with $P_{j,k}\left(\big|\dejak{j}{v}{k}(t)\big| \right)^2$ is the graph of a map $f:[0,(v^j_j)^2]\to \R_{>0}$ (see the proof of (1) in Proposition \ref{propo: props curva vjk}). An example of this situation can be seen in Figure \ref{fig curva y segmento}. The fact that $2 \cos(t) |v^j_k| \sin(t) \sqrt{(v^k_k)^2-|v^j_k|^2}$ is strictly concave for $t\in(0,\pi/2)$ and that the graph of the map $f$ was obtained adding this coordinate to the second coordinate of the mentioned segment implies the concavity of $f$.
		\begin{figure}
			\begin{center}
				\epsfig{file=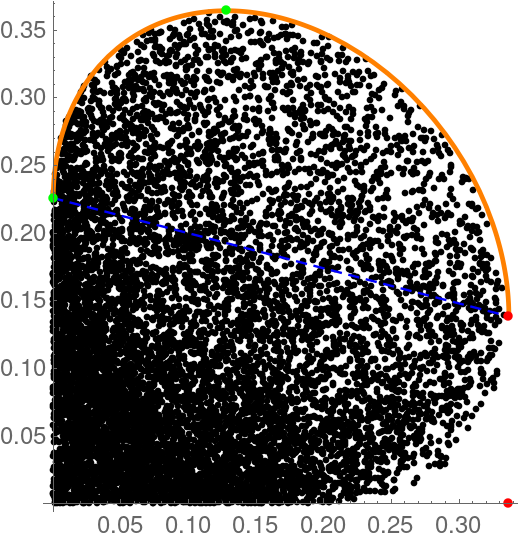, width=5cm}
			\end{center} 
			\caption{Curve $P_{j,k}\left(\big|\dejak{j}{v}{k}(t)\big|^2 \right)$ for $t\in[0,\pi/2]$ in orange, segment between $((v^j_j)^2,|v^j_k|^2)$ and $(0,(v^k_k)^2-|v^j_k|^2)$ dashed in blue, red dots are $((v^j_j)^2,|v^j_k|^2)$ and  $((v^j_j)^2, 0)$, and the green dots are $(0,(v^k_k)^2-|v^j_k|^2)$ and $(|v^k_j|^2,(v^k_k)^2)$. The black dots are projected random points of $m_S$.}
			\label{fig curva y segmento}
		\end{figure}
		Then for $t\in(0,\pi/2)$, if there exists $x,y\in S$, with $\|x\|=\|y\|=1$, and $0\leq \lambda\leq 1$ such that 
		\[
		P_{j,k}\left(\big|\dejak{j}{v}{k}(t)\big| \right)^2=\lambda (|x_j|^2,|x_k|^2)+(1-\lambda) (|y_j|^2,|y_k|^2)
		\]
		then there exist $t_x$ and $t_y$ such that \eqref{eq: igualdad y desig de xj xk en mom} holds for $x$ and $y$. But if any of the inequalities in the $k^\text{th}$ coordinate is strict, then the concavity of the map $f$ is contradicted. Therefore all must be equalities and then $P_{j,k}\left(\big|\dejak{j}{v}{k}(t)\big| \right)^2$ is an extreme point in  $P_{j,k}(m_S)$.
		This in turn implies that $|\dejak{j}{v}{k}(t)|^2$ is in the boundary of $m_S$ for $t\in(0,\pi/2)$.
		
	\item Case $v^j\perp v^k$
		
		As we have seen in \eqref{eq: vjk cero sii vj perp vk} this is equivalent to $v^j_k=0$.
		
		Then 
		$
		P_{j,k}\left(\big|\dejak{j}{v}{k}(t)\big| \right)^2 =\cos^2(t) \left( (v^j_j)^2   \ ,\  |v^j_k|^2\right)+ \sin^2(t)\left(0, \left((v^k_k)^2-|v^j_k|^2\right)\right) 
		$
		which is a segment in $P_{j,k}(m_S)$. If there exists $x,y\in S$, with $\|x\|=\|y\|=1$, and $0\leq \lambda\leq 1$ such that 
		\[
		P_{j,k}\left(\big|\dejak{j}{v}{k}(t)\big| \right)^2=\lambda (|x_j|^2,|x_k|^2)+(1-\lambda) (|y_j|^2,|y_k|^2).
		\]
		we can choose $t_x, t_y$ as in the previous case such that \eqref{eq: igualdad y desig de xj xk en mom} holds. Now if any of the two inequalities given by  \eqref{eq: igualdad y desig de xj xk en mom} is strict there must be $t_0$ such that $P_{j,k}\left(\big|\dejak{j}{v}{k}(t_0)\big| \right)^2$ is not in the segment, which is a contradiction. Then, all inequalities of \eqref{eq: igualdad y desig de xj xk en mom}  are equalities and therefore $P_{j,k}\left(\big|\dejak{j}{v}{k}(t)\big| \right)^2$ is a boundary point of $P_{j,k}(m_S)$. Then $\big|\dejak{j}{v}{k}(t)\big|$ is a boundary point of $m_S$ for $t\in(0,\pi/2)$.
	\end{enumerate}
\end{proof}
\begin{figure}
	\begin{center}
		\epsfig{file=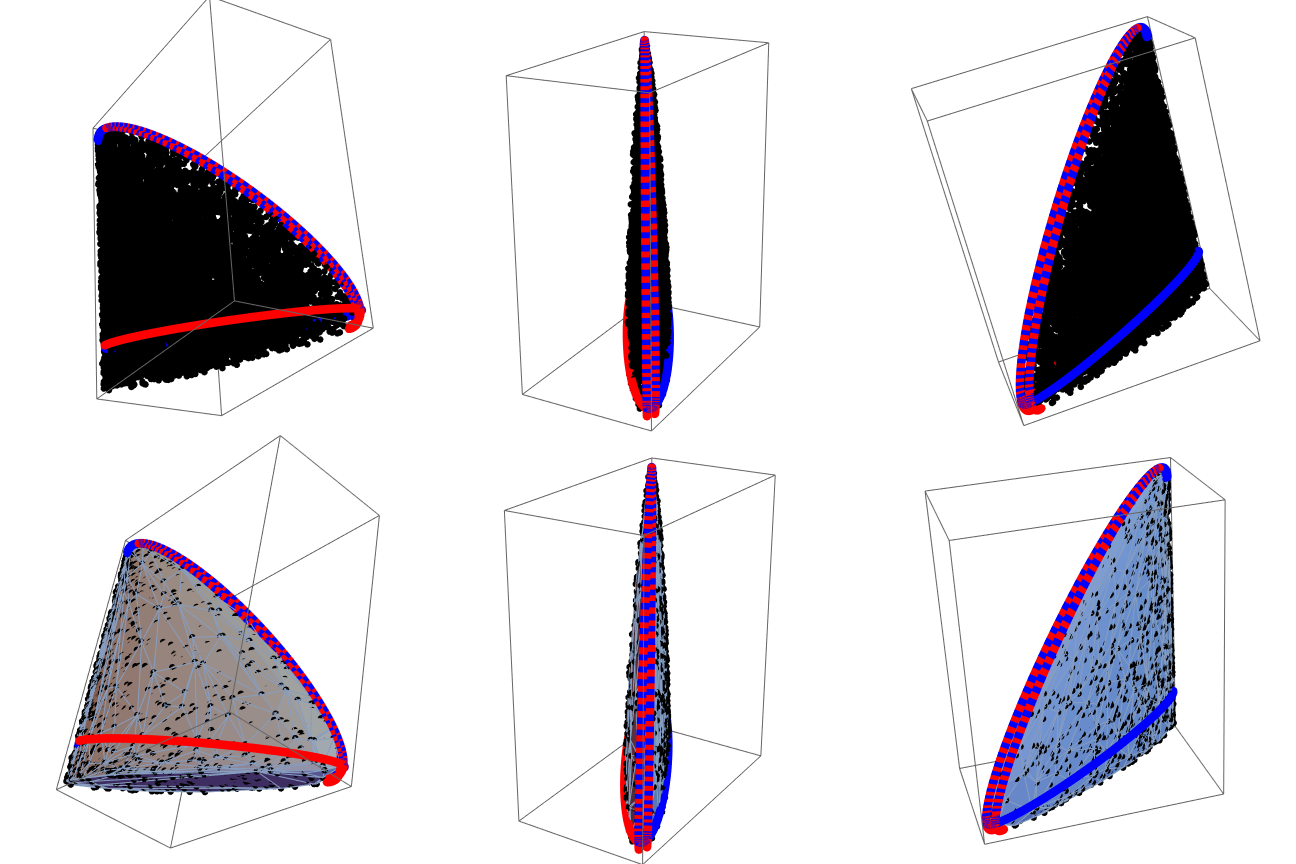, width=16cm	}
	\end{center} 
			\caption{For a subspace $S\subset \C^7$ with dim$(S)=3$, this is an approximation (using random points above and its convex hull below) of the	
				projection of $m_S$ to $\R^3_{\geq 0}$ in 3 fixed coordinates and the plot of 6 curves $\dejak{j}{v}{k}$ (see \eqref{eq: def vtildejk} coinciding with the $j, k$ coordinates being projected. Note the overlapping of $\left|\dejak{j}{v}{k}\right|^2$ and $\left|\dejak{k}{v}{j}\right|^2$ in some parts (see Theorem \ref{vjmoniok a la 2da es igual a vkmonioj a la 2da}).}
			\label{fig curvas estremales en momento proyectado}
\end{figure}
\begin{theorem} 
%
	Let $S$ be a generic subspace of $\C^n$, $v^j$, $v^k$ linearly independent principal vectors such that $v^j \notperp v^k$, $m_S$ the moment of $S$ and the curve $\dejak{j}{v}{k}:[0,\pi/2]\to S$ from Definition \ref{def: jflechakV}.
	
	Then the elements $|\dejak{j}{v}{k}(t)|^2$ are extremal points of $m_S$ for every $ t \in [0,\pi /2)$.
	
	Moreover, if $s|w|^2+(1-s)|z|^2= |\dejak{j}{v}{k}(t_0)|^2$ with $0 < s < 1$, $\|w\|=\|z\|=1$, $w, z\in S$, $ t_0\in [0,\pi /2)$, then $w$ and $z$ must be multiples of $\dejak{j}{v}{k}(t_0)$.
\end{theorem}
\begin{proof}
The case $t_0=0$ has been proved in Proposition \ref{prop: extremalidad de la proyeccion canonica}.

Suppose there exist $w,z \in S$ with $\|w\|=\|z\|=1$ such that  
\begin{equation}
\label{eq: comb. convexa |z|2 y |w|2}
s|w|^2+(1-s)|z|^2=\left|\dejak{j}{v}{k}(t_0)\right|^2\ \text{ with }\ 0 < s < 1
\end{equation}
Then, since $v^j\notperp v^k$ using (2) of Theorem \ref{teo props proy vjk en coords jk} follows that $|z_j|=\left|\dejak{j}{v_j}{k}(t_0)\right|$ and $|z_k|=\left|\dejak{j}{v_k}{k}(t_0)\right|$.

As it was done in \ref{eq: descomp x = a b c} under the notation \ref{notacion wjkMonio}, we can write $z$ from 
\eqref{eq: comb. convexa |z|2 y |w|2} as a linear combination of $v^j,\, \widetilde{w}^{jk}$ and a vector $y\in S$ with $\|y\|=1$, orthogonal to the subspace span$\{v^j, \widetilde{w}^{jk}\}$
\[
z = a\ v^j + b\ \widetilde{w}^{jk} + c\ y \text{ , with } a, b , c \in \C.
\]
Recall that as in the proof of Proposition \ref{prop: desigs entre coords vj xj vk xk} follows that $\widetilde{w}_j^{jk}=y_j=0$ and $y_k=0$, and then it must be $|a|=\cos(t_0)$.
Now observe that
\begin{equation}
\label{eq: igualdad mod zk con mod vjkk}
\begin{split}
|z_k|&=\left| \cos(t_0) e^{i\,\arg(a)}\, v^j_k+b\,  \widetilde{w}_k^{jk}\right|
\leq 
\cos(t_0) \, |v^j_k| + |b|\,  \left|\widetilde{w}_k^{jk}\right|
\\
&\leq
\cos(t_0) \, |v^j_k| + \sin(t_0)\,  \left|\widetilde{w}_k^{jk}\right|
=
\left| \dejak{j}{v_k}{k}(t_0)\right|
\end{split}
\end{equation}
where in the last inequality we used that $|z|^2=\cos^2(t_0)+|b|^2+|c|^2=1$ which implies $|b|\leq \sin(t_0)$.

Since we are supposing \eqref{eq: comb. convexa |z|2 y |w|2} the extremes of the inequality obtained in \eqref{eq: igualdad mod zk con mod vjkk} are equal  and then the complex numbers $\cos(t_0)e^{i\,\arg(a)}\,v^j_k$ and $b\, \widetilde{w}_k^{jk}$ must be collinear. This last statement can only be true if $\arg(a)=\arg(b)$. Moreover, since the inequalities in \eqref{eq: igualdad mod zk con mod vjkk} are equalities, follows that $|b|=\sin(t_0)$, and then $c=0$, which implies that
$
a= e^{i \arg(a)} \cos(t_0) \text{ and } b =e^{i \arg(a)} \sin(t_0)
$. 
Then
\begin{equation}
\begin{split}
z&= a\, v^j+b\, \widetilde{w}^{jk}
= e^{i \arg(a)} \left(\cos(t_0)  v^j + \sin(t_0)\,
\,e^{i \arg\left(v^j_k\right)} \frac{\left(v_k^k -\langle v^k,v^j\rangle v_k^j\right)}{\|v^k -\langle v^k,v^j\rangle v^j\|}\right)\\
&=e^{i\,\arg(a)}\dejak{j}{v}{k}(t_0)
\end{split}
\end{equation}

This proves in particular that $|z|^2=|\dejak{j}{v}{k}(t_0)|^2$. 

Following the same steps it can be proved that $w$ is a multiple of $\dejak{j}{v}{k}(t_0)$, and that $|w|^2= |\dejak{j}{v}{k}(t_0) |^2$, which implies that $ |\dejak{j}{v}{k}(t_0) |^2$ is an extremal point in $m_S$.
\end{proof}


\begin{theorem}\label{vjmoniok a la 2da es igual a vkmonioj a la 2da} Let $S$ be a generic subspace,  $1\leq j, k \leq n $ with $j\neq k$, and $\overset{j \curvearrowright k}{v}$ the curve from Definition \ref{def: jflechakV} but with domain in the interval $ \left[0,\arccos\left( {\left|v_j^k\right|}/{v_j^j}\right)\right]$.
Then 	
	\[
	\text{Dom} (\overset{j \curvearrowright k}{v} )=
	\text{Dom} (\overset{k \curvearrowright j}{v}  )
	\quad \text{ and }\quad
	\text{Im}(| \overset{j \curvearrowright k}{v} |^2 )= \text{Im}( | \overset{k \curvearrowright j}{v} |^2  )
	\]
	(where Dom and Im denote the domain and image of the corresponding curves in $\R^n_{\geq 0}$), and
	for every $t \in \left[0, \arccos({|v_j^k|}/{v_j^j})\right]$ there exists a unique $s=\left(\arccos({|v_j^k|}/{v_j^j})-t\right) \in \left[0, \arccos({|v_k^j|}/{v_k^k})\right]$ such that
	\begin{equation}
	\label{eq: reparametrizacion de t a s}
	\overset{j \curvearrowright k}{v}(t) = {e}^{\text{i} \arg(v_k^j)} \,  \overset{k \curvearrowright j}{v}\left(s\right)
	\end{equation}
	(see Figure \ref{fig curvas estremales en momento proyectado}).
\end{theorem}
	
	\begin{proof}
		The domains of $\overset{k \curvearrowright j}{v}$ and $\overset{j \curvearrowright k}{v} $ are equal since using (\ref{eq: igualdad fracc vjk/vjj etc}) follows that $\frac{| v_k^j |}{v_k^k} = \frac{| v_j^k |}{v_j^j}$.
		
		First observe that $ \overset{j \curvearrowright k}{v}(0) = v^j = {e}^{\text{i} \arg(v_k^j)}\, \overset{k \curvearrowright j}{v}\left( \arccos\left({|v_k^j|}/{v_k^k}\right) \right)$ 
		and that
		$ \overset{j \curvearrowright k}{v}\left( \arccos\left({|v_j^k|}/{v_j^j}\right) \right) = {e}^{\text{i} \arg(v_k^j)}\,v^k = {e}^{\text{i} \arg(v_k^j)}\, \overset{k \curvearrowright j}{v}( 0 )$.
		
		Moreover, the following statements can be proved directly from the properties of $v^j$, $v^k$ and $\dejak{j}{v}{k}$.
		\begin{itemize}
			\item 		$\arg(v^j_k)=\arg\left(\overline{v^k_j}\right)$
			\item $\overset{j \curvearrowright k}{v_j}(t) \in \R $, for all $t \in \left[0,\arccos\left({|v_j^k|}/{v_j^j}\right)\right]$
			\item $\arg\left(\overset{j \curvearrowright k}{v_k}(t)\right) = \arg(v_k^j)$,  
			for all $t \in \left[0,\arccos\left({|v_j^k|}/{v_j^j}\right)\right]$
			\item ${e}^{\text{i} \arg(v_k^j)}\, \overset{k \curvearrowright j}{v_j}(s) \in \R $, for all $s \in \left[0,\arccos\left({|v_k^j|}/{v_k^k}\right)\right]$
 			\item
			$\arg\left(\overset{k \curvearrowright j}{v_j}(s)\right) = -\arg(v_k^j)$, for all $s \in \left[0,\arccos\left({|v_k^j|}/{v_k^k}\right)\right]$
			\end{itemize}
		Then for every $t $ and $s$ in the interval $\left[0,\arccos\left( {\left|v_j^k\right|}/{v_j^j}\right)\right]=\left[0,\arccos\left( {\left|v_k^j\right|}/{v_k^k}\right)\right]$ the following arguments coincide	 
		\begin{equation}\label{coinciden los argumentos}
		\arg\left( \overset{j \curvearrowright k}{v_j}(t) \right) = \arg\left( {e}^{\text{i} \arg(v_k^j)}\, \overset{k \curvearrowright j}{v_j}(s) \right) \quad \text{ and } \quad \arg\left(\overset{j \curvearrowright k}{v_k}(t)\right) =\arg\left( {e}^{\text{i} \arg(v_k^j)}\, \overset{k \curvearrowright j}{v_k}(s)\right). 
		\end{equation} 
		The real map $\overset{j \curvearrowright k}{v_j}(t)=\cos(t)v_j^j$ is decreasing for $t \in \left[0,\arccos({|v_k^j|}/{v_j^j})\right]$ and using \eqref{eq: props curva vjk} we can write
		\[
		\left|\overset{j \curvearrowright k}{v_k}(t) \right|=\left| \cos(t) v^j_k+\sin(t) e^{i \arg\left(v^j_k\right)} \sqrt{(v^k_k)^2-|v^j_k|^2}\right| =  \cos(t) \left| v^j_k\right|+\sin(t) \sqrt{(v^k_k)^2-|v^j_k|^2}.
		\] 
		This map has only one critical point (maximum) in $[0,\frac{\pi}{2}]$ when $t=\arccos\left({|v_j^k|}/{v_j^j}\right)=\arccos\left({|v_k^j|}/{v_k^k}
		\right)$, and therefore is increasing in $\left[0,\arccos({|v_k^j|}/{v_k^k})\right]$. 
		
		Similarly $\overset{k \curvearrowright j}{v_k}(s)$ is decreasing in $\left[0,\arccos({|v_k^j|}/{v_k^k})\right]$, while  $\left|\overset{k \curvearrowright j}{v_j}(s)\right|$ is increasing in such interval (see for example Figure \ref{fig: curva proyectada solo con modulos}).
 		
		Then if $0\leq s_0 < s_1 \leq \arccos\left({|v_k^j|}/{v_k^k}\right)=\arccos\left({|v_j^k|}/{v_j^j}\right)$ 
		\begin{equation}\label{crec decrec de j y de k}
		\left|\overset{k \curvearrowright j}{v_j} (s_0)\right| <\left|\overset{k \curvearrowright j}{v_j} (s_1)\right| \quad \text{ and } \quad \overset{k \curvearrowright j}{v_k} (s_0) > \overset{k \curvearrowright j}{v_k} (s_1) \quad\text{ hold.}
		\end{equation} 
		  Now given a fixed $s_0 \in [0,\arccos({|v_k^j|}/{v_k^k})]$ if we apply Proposition \ref{prop: desigs entre coords vj xj vk xk} to $x=\overset{k \curvearrowright j}{v} (s_0) \in S$ we can find $ t_0 \in [0,\arccos({|v_j^k|}/{v_j^j})]$ such that
		\begin{equation}\label{s0 menor a t0}
		| \overset{k \curvearrowright j}{v_j} (s_0) | \leq | \overset{j \curvearrowright k}{v_j} (t_0)| \quad \text{and} \quad  | \overset{k \curvearrowright j}{v_k} (s_0) | \leq | \overset{j \curvearrowright k}{v_k} (t_0)|.
		\end{equation}
		Then applying again Proposition \ref{prop: desigs entre coords vj xj vk xk} for the vector $\overset{j \curvearrowright k}{v} (t_0) \in S$, we can find $ s_1 \in [0,\arccos({|v_k^j|}/{v_k^k})]$ such that
		\begin{equation}\label{t0 menor a s1}	   
		| \overset{j \curvearrowright k}{v_j} (t_0)| \leq | \overset{k \curvearrowright j}{v_j} (s_1) | \quad \text{and} \quad  | \overset{j \curvearrowright k}{v_k} (t_0)| \leq | \overset{k \curvearrowright j}{v_k} (s_1) | 
		\end{equation}
		
		Now considering \eqref{crec decrec de j y de k}, equations \ref{s0 menor a t0} and \ref{t0 menor a s1} imply that $s_0=s_1$, and then from \ref{coinciden los argumentos} we obtain that
		\[ 
	    {e}^{\text{i} \arg(v_k^j)}\, \overset{k \curvearrowright j}{v_j} (s_0) = \overset{j \curvearrowright k}{v_j} (t_0) 
		\quad \text{and} \quad 
		{e}^{\text{i} \arg(v_k^j)}\,   \overset{k \curvearrowright j}{v_k} (s_0) = \overset{j \curvearrowright k}{v_k} (t_0).
		\]
			
		The vectors $ {e}^{\text{i} \arg(v_k^j)}\overset{k \curvearrowright j}{v}(s_0)$ and $\overset{j \curvearrowright k}{v} (t_0)$ are equal in their $j^\text{th}$ and $k^\text{th}$ coordinates in a dimension 2 subspace spanned by the vectors $v^j$ and $v^k$ $\in S$.
		Then, for each $z\in\text{gen}\{v^j,v^k\}$, if we consider the linear system $\alpha v^j+\beta e^{i \arg(v^j_k)} v^k = \left(\begin{smallmatrix}
		v^j\\e^{i \arg(v^j_k)} v^k
		\end{smallmatrix}\right)^t\cdot 
		\left(\begin{smallmatrix}
		\alpha\\ \beta
		\end{smallmatrix}\right)=z$, we know that it has a unique solution $\left(\begin{smallmatrix}
		\alpha\\ \beta
		\end{smallmatrix}\right)\in \C^2$ since $v^j$ and $e^{i \arg(v^j_k)} v^k$ are linearly independent.
		Considering only the $j^\text{th}$ and $k^\text{th}$ coordinates of the linear system we can conclude that if $
		\det\left(\begin{smallmatrix} v^j_j & v^k_j\\ v^j_k & v^k_k\end{smallmatrix}\right)\neq 0$, then
		$
		  {e}^{\text{i}\, \arg(v_k^j)}\overset{k \curvearrowright j}{v}(s_0) = \overset{j \curvearrowright k}{v} (t_0).
		$
		But the case $\det\left(\begin{smallmatrix}v^j_j & v^k_j\\ v^j_k & v^k_k\end{smallmatrix}\right)= 0$ is never possible. Otherwise  $ v^j_j v^k_k- v^j_k v^k_j = v^j_j v^k_k- |v^j_k| \ |v^k_j| = 0$ (recall that $\arg(v^j_k)=-\arg(v^k_j)$). Then using Lemma \ref{lem: extremalidad de la coordenada} we obtain that $v^j_j = |v^k_j|$ and $v^k_k=|v^j_k|$, and then that $v^j=e^{i \arg(v^j_k)} v^k$, a contradiction.		
				Therefore for every $s_0\in [0,\arccos({|v_k^j|}/{v_k^k})]$ there exists a unique  $t_0\in [0,\arccos({|v_j^k|}/{v_j^j})]$ such that
				\begin{equation}
				\label{eq vkj es mult de vjk}
				{e}^{\text{i} \, \arg(v_k^j)}\overset{k \curvearrowright j}{v}(s_0) = \overset{j \curvearrowright k}{v} (t_0).
				\end{equation}

To prove the formula \eqref{eq: reparametrizacion de t a s} recall that the domains of $\overset{k \curvearrowright j}{v}$ and $\overset{j \curvearrowright k}{v} $ are equal.
Now take $t\in [0, \arccos({|v_j^k|}/{v_j^j})]$ and define $s = \arccos\left({|v_j^k|}/{v_j^j}\right)-t$. Then
\begin{equation}
\begin{split}
\cos(t)&= \cos(\arccos({|v_j^k|}/{v_j^j}) - s) = \cos(\arccos({|v_j^k|}/{v_j^j})) \cos(s)+ \sin(\arccos({|v_j^k|}/{v_j^j}) ) \sin( s)
\\
&=\left({|v_j^k|}/{v_j^j}\right)  \cos(s)+ \sqrt{1-\left({|v_j^k|}/{v_j^j} \right)^2}  \sin( s)=
\frac{|v_j^k|   \cos(s)+ \sqrt{(v_j^j)^2-|v_j^k|^2}   \sin( s)}{v_j^j}
\end{split}
\end{equation}
Therefore
\begin{equation}
\begin{split}
\overset{j \curvearrowright k}{v_j}(t)&=
\cos(t)  v_j^j = |v_j^k|  \cos(s)+ \sqrt{(v_j^j)^2-|v_j^k|^2} \, \sin( s)
\\
&= \left|v_j^k \, \cos(s)+ \sqrt{(v_j^j)^2-|v_j^k|^2} \, e^{\text{i} \arg(v_j^k)} \sin( s)\right|= e^{\text{i} \arg(v_k^j)}
\overset{k \curvearrowright j}{v_j}(s) 
\end{split}
\end{equation}

The proof for the $k^\text{th}$ coordinate is similar. 

As was seen previously in the proof of \eqref{eq vkj es mult de vjk} this equality in the $j$ and $k$ coordinates is enough to prove the formula \eqref{eq: reparametrizacion de t a s}.
\end{proof}

\section{Relation of $m_S$ with the joint algebraic numerical range}\label{secc JNR}

We will consider here the relation of the moment $m_{S}=m_{S,E}$ (see Definition \ref{def: momento de un subespacio}) of a subspace $S$ related to a fixed orthonormal basis $E$ with the \textsl{joint numerical range}\cite{plaumann-sinn-weis,bonsall-duncan,weis2017}). 
The joint numerical range of $m$ hermitian matrices $A_1, A_2,\dots,A_m\in M_n^h(\C)$, sometimes called the joint algebraic numerical range \cite{muller}, is defined by
\begin{equation}
\label{def JNR}
W(A_1, A_2,\dots,A_m)=
\{\left(\tr(A_1\rho), \tr(A_2\rho),\dots,\tr(A_m\rho)\right)\in\R^n: \rho \in \DD \}
\end{equation}
with $\DD=\{\rho \in M_n^h(\C): \rho\geq 0, \tr(\rho)=1\}$ the set of density matrices of $M_n(\C)$ and $\tr$ the usual trace (sum of diagonal entries).
This set is the convex hull of the also called joint numerical range in the literature (that we will denote classic joint numerical range)
\begin{equation}
\label{def JNR clasico}
\begin{split}
W_\text{class}(A_1, A_2,\dots,A_m)&=\\
=\{\left(\langle A_1 x,x\rangle,\right.&\left.\hspace{-.15cm}\langle A_2 x,x\rangle,\dots,\langle A_m x,x\rangle\right)\in\R^n: x \in \C^n, \|x\|=1 \}
\end{split}
\end{equation}
that is not necessarily convex in general. In the case $W_{\text{class}}$ is convex then the equality $W=W_{\text{class}}$ holds.

More precisely, we will see that $m_{S}$ can be described as a particular subset of a joint numerical range of some selected hermitian matrices (see Theorem \ref{teo mV es JNR y suma 1}).

Recall that one of the equivalent definitions of $m_S$ seen in \eqref{defequiv mS 1 en prop} of Proposition \ref{prop: equivalencias de momento} is considering the following identification of the vectors of $m_S$ with the diagonal entries of certain matrices
\[
m_{S,E}  \simeq \text{Diag}\{Y \in M_n^h(\C):Y\geq 0, \tr(Y)=1, \text{Im}(Y) \subset S\}.
\] 
Since these $Y$ are hermitian matrices the condition Im$(Y)\subset S$ is equivalent to $Y=P_S Y=Y P_S=P_S Y P_S$. This allow us to rewrite the moment set $m_S$ as
\begin{equation}\label{eq donde se relaciona mV con diags de Ys}
\mS\simeq\text{Diag}\left(\{Y\in M_n^h(\C):Y\in\DS\}\right)
\end{equation}
where $\DS$  
\begin{equation}\label{def DsubV}
\begin{split}
\DS&= \left\{Y\in M_n^h(\C):Y\geq 0, \tr(Y)=1, P_S Y=Y(=YP_S)\right\}\\
&=\left\{Y\in \DD: Y=P_S Y=YP_S=P_S  Y P_S\right\}.
\end{split}
\end{equation}

	\begin{lemma}
	Consider $\DS$ as in \eqref{def DsubV} and ${0}\subsetneq S\subsetneq \C^n$ a subspace of $\C^n$. Then  
	\[
	\DS= \left\{\rho\in\DD: \tr(P_S\rho P_S)=1\right\}.
	\]
\end{lemma}
\begin{proof}
	Observe that 
	
	\begin{equation*}
	\begin{split}
	\rho\in\DD_S &\Leftrightarrow 
	\rho\geq 0 , \tr(\rho)=1, \rho=P_S \rho=\rho P_S=P_S \rho P_S\\
	&\Leftrightarrow \rho\geq 0 , \tr(\rho)=1, \rho =\frac{1}{\tr(P_S \rho P_S)}P_S \rho P_S\\
	&\Leftrightarrow \rho\geq 0 , \tr(P_S \rho P_S)\neq 0, \rho =\frac{1}{\tr(P_S \rho P_S)}P_S \rho P_S	\end{split}
	\end{equation*}
	and then
	\begin{equation}\label{eq: equiv def de DsubV} 
	\DS=\left\{\rho\in\DD: \tr(P_S\rho P_S)\neq 0, \rho=\frac1{\tr(P_S\rho P_S)}{P_S\rho P_S}\right\}
	\end{equation}
	which implies that $\DS\subset \left\{\rho\in\DD: \tr(P_S\rho P_S)=1\right\}$.
	
	Now consider a $\rho\in\DD$ with a spectral decomposition $\rho=\sum_{i=1}^k p_i\ x_i\cdot (x_i)^*$ with $\sum_{i=1}^k p_i=1$, $p_i>0$ and $x_i\in\C^n$, $\|x_i\|=1$, $x_i\perp x_j$ (if $i\neq j$). Then the property $\tr(P_S \rho P_S)=1$ implies that $1=\sum_{i=1}^k p_i \ \tr\left(P_S x_i\cdot (x_i)^*P_S\right)=\sum_{i=1}^k p_i\ \tr\left(P_S x_i\cdot \left(P_S x_i\right)^*\right)$. Then it must be $\tr(P_S x_i\cdot \left(P_S x_i\right)^*)=1$ for all $i=1,\dots,k$ which in turn implies that $\|P_S x_i\|=1$ and therefore $x_i\in S$ and  $P_S\,\rho=\rho=\rho^*=\rho\, P_S=P_S\, \rho\, P_S$. We have obtained that
	\begin{equation}
	\label{eq: DV en funcion de funciones de densidad}
	\left\{\rho\in\DD: \tr(P_S\rho P_S)=1\right\}\subset \DS.
	\end{equation}
	which concludes the proof.
\end{proof}

\begin{notation} \label{notation 2} Consider a fixed basis $E=\{e_1, e_2,\dots,e_n\}$ of $\C^n$.
	The rank one orthogonal projections onto the subspace generated by a single $e_i\in E$ are described in the same $E$ basis by the $n\times n$ matrices $E_i=e_i\cdot (e_i)^*$. Here $\cdot$ denotes the matrix product of the column vector $e_i$ with the row vector $(e_i)^*$(conjugate transpose of $e_i$). In this case, the coordinates of the vector $e_i$ are zeros with the exception of a $1$ in its $i^\text{th}$ coordinate and $E_i$ is a $n\times n$ matrix of zeros and only a $1$ in its $i, i$ entry. 
	In this case, the $E_i$ projections are also denoted with $e_i\otimes e_i$. 
\end{notation}

\begin{lemma}\label{lema: momento igual a trazas EiY}
	If $\{0\}\subsetneq S\subsetneq \C^n$ is a subspace of $\C^n$ and $\DS$ is as in \eqref{def DsubV}, then
	\[
	\mS=\{\left(\tr(E_1 Y), \tr(E_2 Y),\dots,\tr(E_n Y)\right)\in\R^n_{\geq 0}: Y \in \DS \}
	\]
	where we denote with $E_i= e_i\otimes e_i$ (see Notation \ref{notation 2}) and $\tr$ is the usual trace.
\end{lemma}
\begin{proof}
	Let $Y\in\DS$ and a generic element $\text{Diag}(Y)\in m_S$ (see \eqref{eq donde se relaciona mV con diags de Ys} and \eqref{def DsubV}) under the identification $\simeq$ between diagonals of positive definite matrices and vectors in $\R^n_{\geq 0}$. Then, if  $E_i=e_i\otimes e_i$, the proof follows observing that 
		\begin{equation*}
		\begin{split}
		m_S\ni\Diag(Y)&\simeq (\tr(E_1 Y E_1), \tr(E_2 Y E_2),\dots,\tr(E_n Y  E_n))
		\\
		&= (\tr(E_1 Y), \tr(E_2 Y),\dots,\tr(E_n Y))
			\end{split}
		\end{equation*} 
		\end{proof}

The next result shows the relation between $m_S$ and a particular joint numerical range.

\begin{theorem}\label{teo mV es JNR y suma 1}
	Let  $E=\{e_i\}_{i=1}^n$ be a fixed orthonormal basis, $m_S$ the moment of a subspace $S$ with $\{0\}\subsetneq S\subsetneq\C^n$ related to $E$ (see \eqref{def: momento de un subespacio}), $P_S$ the orthogonal projection onto $S$, $E_i=e_i\otimes e_i$ the rank one orthogonal projections onto the subspace generated by $e_i\in E$ (see details in Notation \ref{notation 2}) and $W(A_1,A_2,\dots,A_n)$ the joint numerical range of the hermitian matrices $A_1,A_2,\dots,A_n\in M_n^h(\C)$ (see \eqref{def JNR}). Then
	\begin{equation}\label{eq: mV es JNR interseccion suma 1}
	\mS=W\left(P_S E_1 P_S , P_S E_2 P_S ,\dots, P_S E_n P_S\right)\ \bigcap\ \left\{x\in \R^n:\sum_i x_i=1\right\}.
	\end{equation}
\end{theorem}
\begin{proof}
    Suppose $x\in m_S$. Then using Lemma \ref{lema: momento igual a trazas EiY} we can find $Y\in\DS$ such that
	$x=(\tr(E_1 Y) , \tr(E_2 Y), \dots, \tr(E_nY))$. Then if we consider that $Y=YP_S=P_S Y=YP_S Y$ follows that
	\begin{equation}\label{eq mom en jnr}
	\begin{split}
	x&=\left(\tr(E_1P_S Y P_S),\tr(E_2P_S Y P_S),\dots, \tr(E_1P_S Y P_S)\right)\\
	&=\left(\tr(P_S E_1P_S Y ),\tr(P_S E_2P_S Y),\dots, \tr(P_S E_1P_S Y)\right).
	\end{split}
	\end{equation}
	Then $Y\in \DS\subset \DD$ and \eqref{eq mom en jnr} imply that $x\in W\left(P_S E_1P_S ,P_S E_2P_S,\dots, P_S E_1P_S \right)$. Moreover, if $x=(x_1,x_2,\dots,x_n)$ then
	\[
	\sum_{i=1}^n x_i=\sum_{i=1}^n \tr(E_iY)=\tr\left(\sum_{i=1}^n E_iY\right)=\tr\left(Y\sum_{i=1}^n E_i\right)=
	\tr\left(Y\, I\right)=\tr(Y)=1
\]
	which proves that $m_S$ is included in $W\left(P_S E_1P_S ,P_S E_2P_S,\dots, P_S E_1P_S \right)$ and its coordinates add to one.
	
	Let us prove the other inclusion. If $x\in W\left(P_S E_1P_S ,P_S E_2P_S,\dots, P_S E_1P_S \right)$ and its coordinates add to one, then there exists $\rho\in \DD$ such that
	\begin{equation}
	\begin{split}
	x&=\left(\tr(P_S E_1 P_S \rho) , \tr(P_S E_2 P_S \rho),\dots,\tr(P_S E_n P_S \rho)\right)\\
	&=\left(\tr(E_1 P_S \rho P_S) , \tr(E_2 P_S \rho P_S),\dots,\tr(E_n P_S \rho P_S)\right)\\
	\end{split}
	\end{equation}
	and $1=\sum_{i=1}^n \tr(E_i P_S \rho P_S) =\tr(P_S \rho P_S\sum_{i=1}^n E_i) =\tr(P_S \rho P_S)$. 
	Then if we consider $Y=P_S \rho P_S$ it is apparent that $Y\in\DS$ and therefore by Lemma \ref{lema: momento igual a trazas EiY} follows that $x= (\tr(E_1 Y) , \tr(E_2 Y) , \dots , \tr(E_nY))\in m_S$.
\end{proof}

The following result is apparent after considering that $W_\text{class}=W$ if and only if $W_\text{class}$ is convex.
\begin{corollary}\label{coro Wclass convexo y mS}
	Let us suppose that $S$ satisfies the assumptions of the previous Theorem \ref{teo mV es JNR y suma 1}. 
	
	Then if $W_{\text{class}}\left(P_S E_1 P_S , P_S E_2 P_S ,\dots, P_S E_n P_S\right)$ is convex,
	\[
	\mS=W_\text{class} \left(P_S E_1 P_S , P_S E_2 P_S ,\dots, P_S E_n P_S\right)\ \bigcap\ \left\{x\in \R^n:\sum_i x_i=1\right\}
	\]
	and 
	\[
	\text{cone}\left(m_S\right) =\text{cone}\left(W_\text{class}\left(P_S E_1 P_S , P_S E_2 P_S ,\dots, P_S E_n P_S\right)\right).
	\]
\end{corollary}

Moreover, $m_S$ and $W\left(P_S E_1 P_S , P_S E_2 P_S ,\dots, P_S E_n P_S\right)$ generate the same positive cone as showed in the following lemma (see \eqref{eq: def cono conj} for a definition of positive cone).
\begin{proposition}\label{prop cono igual cono}
	Let $\{0\}\subsetneq S \subsetneq \C^n$ be a subspace of $\C^n$, then
	\[
	\text{cone}\left(m_S\right) =\text{cone}\left(W\left(P_S E_1 P_S , P_S E_2 P_S ,\dots, P_S E_n P_S\right)\right)
	\]
	where we use the definition of positive cone generated by a set given in \eqref{eq: def cono conj}.
\end{proposition}	
\begin{proof}
The inclusion 
\[
\text{cone}\left(m_S\right)\subset \text{cone}\left(W\left(P_S E_1 P_S ,  \dots, P_S E_n P_S\right)\right)
\]
 is obvious since using Theorem \ref{teo mV es JNR y suma 1} follows that $m_S\subset W\left(P_S E_1 P_S ,  \dots, P_S E_n P_S\right)$.

Note that the null vector belongs to both cones. Suppose that $c\in \text{cone}\left(W\left(P_S E_1 P_S ,  \dots, P_S E_n P_S\right)\right)$ and $c\neq (0,0,\dots,0)$. Then there exists $t>0$, $\rho\in\DD$ such that
\[
c=t \left(\tr(P_S E_1 P_S \rho) , \tr(P_S E_2 P_S \rho),\dots,\tr(P_S E_n P_S \rho)\right)
\]
Now since $c\neq (0,0,\dots,0)$ and $c_i=\tr(P_S E_i P_S \rho)=\tr( \rho^{1/2}P_S E_i P_S \rho^{1/2})\geq 0$ then $\sum_{i=1}^n\tr(P_S E_i P_S \rho)=\tr\left( P_S \rho P_S \sum_{i=1}^n E_i\right)=\tr( P_S \rho P_S)>0$. Let us denote with $\tau=\tr( P_S \rho P_S)$ and observe that
\begin{equation}
\label{eq c multiplo posit de mS}
c= t\ \tau \left(\tr\left( E_1 (1/\tau) P_S \rho P_S\right) , \tr\left( E_2 (1/\tau) P_S \rho P_S\right),\dots,\tr\left( E_n (1/\tau) P_S \rho P_S\right)\right)
\end{equation}
with $(1/\tau) P_S \rho P_S\geq 0$ and $\tr\left((1/\tau) P_S \rho P_S\right)=1$. This implies that $(1/\tau) P_S \rho P_S\in \DS$ and therefore $c$ is a positive multiple of an element of $m_S$ (see \eqref{eq c multiplo posit de mS} and Lemma \ref{lema: momento igual a trazas EiY}).
\end{proof}
The following result gives a bound of how close are two moments sets of subspaces when their respective projections are close in norm.
\begin{proposition}
Consider a pair of subspaces $S$ and $V$ of $\C^n$ and the notations of $\DS$ as in \eqref{def DsubV}, $m_S$, $E$, $E_i$ for $i=1,\dots,n$ and the joint numerical range $W$ as in Theorem \ref{teo mV es JNR y suma 1}. Then, if $\|P_S-P_V\|< \frac1{2n}$, 
\begin{equation}
\begin{split}
\dist_H(m_S,m_V)&=\max \left\{ \sup_{x\in m_S} \dist(x,m_V) , \sup_{z\in m_V} \dist(z,m_S)\right\} \\
&\leq  (2 \sqrt{n}+1) \|P_S-P_V\|
\end{split}\end{equation}
where $\dist_H$ denotes the Hausdorff distance of sets.
\end{proposition}
	\begin{proof}
		We will first prove that if $x\in m_S$ then $\dist(x, m_V)<  (2 \sqrt{n}+1) \|P_S-P_V\|$. 
		
		Using Lemma \ref{lema: momento igual a trazas EiY} if $x\in m_S$ then there exists $Y\in \DS$ ($Y\geq 0$, $\tr(Y)=1$, $P_S Y=Y$) such that  $x=(\tr(E_1 Y),  \dots,\tr(E_n Y))$. 
		
		Now consider the vector $(\tr(E_1 P_V Y P_V),  \dots,\tr(E_n  P_V Y P_V))$ and observe that 
		\[
		\sum_{i=1}^n \tr(E_i P_V Y P_V)=\sum_{i=1}^n (P_V Y P_V)_{i,i}=\tr(P_V Y P_V)\geq 0.
		\]
		But that sum cannot be zero. Otherwise, note that $\tr(P_V Y P_V)=0$ if and only if $P_V Y P_V=0$ (recall that $Y\geq 0$). Then using the assumption that $\|P_S-P_V\|< \frac1{2n}$ follows that
		\[
		\|P_S Y P_S-P_V Y P_V\| \leq \|P_S Y (P_S-P_V)\| + \|(P_S-P_V) Y P_V\|\leq 2 \|P_S-P_V\| <\frac1n.
		\]
        This would imply that $\|P_S Y P_S-0\|=\|Y\|<1/n$ which contradicts the assumptions that $Y\geq 0$ and $\tr(Y)=1$.
		
		Using again that $\tr(Y)=1$
		\begin{equation}
		\begin{split}
		\label{eq tr PVYPV - 1 chico}
		|\tr(P_V Y P_V)-\overbrace{\tr(P_S Y P_S)}^{1}| &\leq |\tr(P_V Y  )-\tr(P_S Y  )| \leq |\tr\left((P_V-P_S) Y  \right)| \\
		&\leq  \|P_S-P_V\|
		\end{split}
		\end{equation}
		which implies that $1- \|P_S-P_V\|\leq \tr(P_V Y P_V)\leq 1+ \|P_S-P_V\|$ and 
		\begin{equation}\label{eq cota dif tr PS y PV}
		\begin{split}
		|\tr(E_i P_V Y P_V)-\tr(E_i P_S Y P_S) | &= 		\left|\tr\left(E_i (P_V Y P_V - P_S Y P_S)\right) \right|
		\\
		&\leq \tr(E_i) \|P_V Y P_V - P_S Y P_S\|
		\leq 2\|P_S-P_V\|.
		\end{split}
		\end{equation}
		
		It is clear that $\frac1{\tr(P_V Y P_V)} P_V Y P_V\in \DV$ and therefore the vector $z= \left(\frac{\tr(E_1 P_V Y P_V)}{\tr(P_V Y P_V)} , \dots,\frac{\tr(E_n P_V Y P_V)}{\tr(P_V Y P_V)}\right)$ belongs to $m_V$ (see Lemma \ref{lema: momento igual a trazas EiY} and Theorem \ref{teo mV es JNR y suma 1}).
		Then
		\begin{equation*}
		\begin{split}
		\|x-z\|  
		&\leq   
		  \| x-\tr(P_V Y P_V) z\| +\|\tr(P_V Y P_V) z-z\| \\
		&\leq   
		\left(\sum_{i=1}^n\left|\tr(E_i(P_S Y P_S-P_V Y P_V))\right|^2 \right)^{1/2} +
		\left| \tr(P_V Y P_V)-1\right| \|z\|\\
		&\leq 2 \sqrt{n} \|P_S-P_V\|+\|P_S-P_V\| =(2 \sqrt{n}+1) \|P_S-P_V\|
		\end{split}
		\end{equation*}
		where in the last inequality we used \eqref{eq cota dif tr PS y PV}, \eqref{eq tr PVYPV - 1 chico} and the fact that $\|z\|\leq 1$ because its coordinates are positive and add to one.
		
		Then we have proved that for every $x\in m_S$ there exists a $z\in m_V$ such that $\|x-z\|\leq (2 \sqrt{n}+1) \|P_S-P_V\|$. Therefore $\dist(x,m_V)\leq (2 \sqrt{n}+1) \|P_S-P_V\|$. We can argue in a similar way to prove that for every $z\in m_V$ holds $\dist(z,m_V)\leq (2 \sqrt{n}+1) \|P_S-P_V\|$. Then we have proved that
		\begin{equation*}
		\begin{split}
		\dist_H(m_S,m_V)&=\max \left\{ \sup_{x\in m_S} \dist(x,m_V) , \sup_{z\in m_V} \dist(z,m_S)\right\} \leq  (2 \sqrt{n}+1) \|P_S-P_V\|
		\end{split}
		\end{equation*}
	\end{proof}
	\begin{remark}
		A reciprocal of the previous proposition does not hold. There exist orthogonal subspaces $S=\text{gen}(x)$ and $V=\text{gen}(\overline{x})$ generated by a vector $x\in \C^n$ and an orthogonal $\overline{x}$ obtained after conjugating the coordinates of $x$ such that $m_S=\{|x|^2\}=m_V$ and $\|P_S-P_V\|=\sqrt{2}$. 
	\end{remark}

The following are basic observations that describe elements of $W\left(P_S E_1 P_S , P_S E_2 P_S ,\dots, P_S E_n P_S\right)$. These can be obtained directly using properties of the trace, the projections and the principal standard vectors.

\begin{remark}\label{rem elems jnr SE}
	Let $S$ be a subspace of $\C^n$,  $\rho\in \DD$, $E_i=e_i\otimes e_i$ for $i=1,\dots,n$, $\{v^i\}_{i=1}^n$ the principal standard vectors and denote with
	\[
	\Delta(\rho)=\left(\tr\left(P_S E_1 P_S\rho\right), \tr\left(P_S E_2 P_S\rho\right), \dots, \tr\left(P_S E_j P_S\rho\right), \dots,\tr\left(P_S E_n P_S\rho\right)\right)
	\]
	the elements of $W\left(P_S E_1 P_S , P_S E_2 P_S ,\dots, P_S E_n P_S\right)$ corresponding to $\rho\in\DD$.
	Then the following statements hold.
	\begin{enumerate}
		\item If $\rho\in\DD$, then
		\[
		\Delta(\rho)=\left( (P_S \rho P_S)_{1,1} , (P_S \rho P_S)_{2,2}, \dots, (P_S \rho P_S)_{n,n}\right).
		\]
		\item If $\rho\in\DD_S$ (see \eqref{def DsubV}), then
		\[
		\Delta(\rho)=\left( \rho_{1,1}  ,  \rho_{2,2} , \dots, \rho_{n,n}\right).
		\]
		\item If $x\in\C^n$ with $\|x\|=1$, then $x\otimes x\in \DD$ and
		\[
		\Delta(x\otimes x)= \left( |(P_S x)_1|^2 , |(P_S x)_2|^2, \dots, |(P_S x)_n|^2\right)
		\]
    	\item 
    	If $x^j\in\C^n$ with $\|x^j\|=1$, for $j=1,\dots, k$, such that $\rho=\sum_{j=1}^k t_j (x^j\otimes x^j)\in \DD$, then
		\[
    	\Delta(\rho)= \sum_{j=1}^k t_j \Delta(x^j\otimes x^j)= \sum_{j=1}^k t_j \left( |(P_S x^j)_1|^2 , |(P_S x^j)_2|^2, \dots, |(P_S x^j)_n|^2\right)
    	\]
    	\item 
    	If $x\in\C^n$ with $\|x\|=1$, $x\otimes x\in \DD$ and if $P_S=U I_r U^*$ is a spectral decomposition of $P_S$ ($r=\dim(S)$, $I_r$ a $n\times n$ matrix with $1$ in the first $r$ entries of its diagonal and zero elsewhere, and $U$ a unitary whose columns are orthonormal eigenvectors of $P_S$, with the first $r$ belonging to $S$), then
	    \begin{equation*}
	    \begin{split}
	    \tr\left(I_r U^* E_i U I_r (x\otimes x)\right) &=
	    \left|\langle (\overline{U_{1,i}},\overline{U_{2,i}},\dots, \overline{U_{r,i}}), (x_1,x_2,\dots,x_r)\rangle\right|^2\\
	    &= \left|  \langle 
	     \text{row}_i\left(U^*\right), (x_1,x_2,\dots,x_r, 0,\dots, 0)\rangle
	    \right|^2 \\
	    &= \left|  \langle 
	    \text{col}_i\left(U\right), (\overline{x_1},\overline{x_2},\dots,\overline{x_r}, 0,\dots, 0)\rangle
	    \right|^2
	    \text{ , for all } i=1,\dots, n.
	    \end{split}
	    \end{equation*}
	    Therefore, the vectors
			    	    \[
			    	    \Delta(x\otimes x)=
			    	    \left|U^* \cdot
			    	    \left(\begin{smallmatrix}
			    	    x_1\\x_2\\ \vdots\\x_r\\ 0\\ \vdots\\ 0
			    	    \end{smallmatrix}
			    	    \right)
			    	    \right|^2
			    	    \]
	    are elements of $W\left(P_S E_1 P_S , \dots, P_S E_n P_S\right)$ after using the unitary invariance of the joint numerical range. Moreover, the convex hull of these elements equals this joint numerical range.
		\item Consider $P_S=U I_r U^*$ as in the previous item, and$\rho=U^*\mu U\in\DD$ for $\mu\in \DD$ (any $\rho\in\DD$ can be written in this format). Then
		\begin{equation*}
		\begin{split}
		\Delta(\rho)&=
		\left(\langle I_r \mu I_r \left(U e_1\right) I_r, I_r(U e_1)I_r\rangle,\dots ,   \dots,
		\langle I_r \mu  I_r \left(U e_n\right) I_r, I_r(U e_n)I_r\rangle \right)\\
		&=
		\left( \left(\begin{smallmatrix}
		\overline{U_{1,1}}&	,\dots,	&\overline{U_{r,1}}
		\end{smallmatrix}	 
			\right)\cdot 
			\left(\begin{smallmatrix}
			\mu_{1,1}&\dots&\mu_{1,r} 
			\\
			\vdots&\ddots&\vdots 
			\\
			\mu_{r,1}&\dots&\mu_{r,r} 
		 	\end{smallmatrix}\right)
			\cdot 
			\left(\begin{smallmatrix}
			U_{1,1}	\\	\vdots	\\	U_{r,1}	
			\end{smallmatrix}\right)
			,\dots, \right.
			\\
			& \hskip3.5cm \dots 
		\left.\left(\begin{smallmatrix}
		\overline{U_{1,n}}&	,\dots,	&\overline{U_{r,n}}
		\end{smallmatrix}	 
		\right)\cdot 
			\left(\begin{smallmatrix}
			\mu_{1,1}&\dots&\mu_{1,r} 
			\\
			\vdots&\ddots&\vdots 
			\\
			\mu_{r,1}&\dots&\mu_{r,r} 
			\end{smallmatrix}\right)
			\cdot 
			\left(\begin{smallmatrix}
			U_{1,n}	\\	\vdots	\\	U_{r,n}	
			\end{smallmatrix}\right)
			\right)
		\end{split}
		\end{equation*}

		\item If $e_i$ is a member of the standard basis $E$ and we consider $e_i\otimes e_i\in\DD$, then
		\[
		\Delta(e_i\otimes e_i)= \left( (v_1^1)^2 |v_i^1|^2 ,\dots, (v_j^j)^2 |v_i^j|^2, \dots, (v_n^n)^2 |v_i^n|^2\right).
		\]
		\item If $s\in S$ with $\|s\|=1$, then $s\otimes s\in\DD_S\subset \DD$ and
		\[
		\Delta(s\otimes s)= \left( |s_1|^2 , |s_2|^2, \dots, |s_n|^2\right).
		\]

		\item If $v^j\in S$ for $j=1,\dots,n$ is a principal standard vector (see Definition \ref{def: principal standard vectors}), then $v^j\otimes v^j\in \DD_S$ and
		\[
		\Delta(v^j\otimes v^j)= \left( |v^j_1|^2 , |v^j_2|^2, \dots, |v^j_n|^2\right).
		\]
	\end{enumerate}
\end{remark}

		\begin{remark}\label{rem props nuestros jnr}
			\begin{enumerate} Let us consider now some properties related to the joint numerical range of these particular matrices. Here $S$ is a dim$(S)=r$ subspace of $\C^n$.
\item \label{item 1 remark 15} Using basic properties of joint numerical ranges 
follows that if $S$ is a generic subspace
\[
W\left(P_S E_1 P_S , ,\dots, P_S E_n P_S\right)=
\left(\begin{smallmatrix}
(v_1^1)^2&0&\dots&0
\\
0&(v_2^2)^2&0&\vdots
\\ 
\vdots&0 &\ddots& \vdots
\\ 
0&\dots &0& (v_n^n)^2
\end{smallmatrix}\right)\cdot \Big( W\left(v^1\otimes v^1, v^2\otimes v^2, \dots,v^n\otimes v^n\right)\Big)
\]
which provides a way to write the first joint numerical range in terms of another involving only rank one projections related to the principal standard vectors of $S$ (see \ref{def: principal standard vectors}).
\item 
Denote with $W_{V,B}=W\left(P_V B_1 P_V , P_V B_2 P_V,\dots, P_V B_n P_V\right)$, for $V$ a subspace of $\C^n$, $B$ an ordered basis of $\C^n$ and $B_i$ the projection onto its $i^\text{th}$ vector. If we consider a unitary matrix $U\in M_n(\C)$, the subspace $U(S)$ and the basis $U \cdot E$ whose elements are the columns of $U$, then
	\begin{equation*}
	\begin{split}
	W_{S,E}&=W\left(P_S E_1 P_S , ,\dots, P_S E_n P_S\right) \\
	&=W\left(U P_S U^* U E_1 U^* U P_S U^*  ,\dots, U P_S U^* U E_n U^* U P_S U^*\right)\\
	&= W_{U(S)\, ,\, U \cdot E}.
	\end{split}
	\end{equation*}
	This follows using that $P_{U(S)}=U P_S U^*$ and the unitary invariance of the joint numerical range $W$.
\item \label{item 3 remark 15}
Consider a unitary matrix $U$ such that $U P_S U^*=I_r$, with $I_r$ the diagonal matrix with $r$ ones at the beginning and zeroes afterwards The first $r$ columns of $U^*$, that we denote $\{s^i\}_{i=1}^r$, form an orthonormal basis of $S$, and the $n-r$ other $o^1,\dots,o^{n-r}$ columns of $U^*$ form an orthonormal basis of $S^\perp$. 
Then $U(S)= \text{gen}\left\{e_i\right\}_{i=1}^r$ and $U\cdot E=B= \{U\cdot e_i\}_{i=1}^n$.
Then using the previous item 
\[
W_{S,E}=W_{\text{gen}\left\{e_i\right\}_{i=1}^r\, ,\, \{U\cdot e_i\}_{i=1}^n}=W(I_r U e_1\otimes I_r U e_1\, ,\,  I_r U e_2\otimes I_r U e_2\, ,\,  \dots,I_r U e_n\otimes I_r U e_n )
\]
which describes $W_{S,E}$ using only the first $r$ coordinates of the vectors $\{U e_i\}_{i=1}^n$.
\end{enumerate}
\end{remark}

The following results describe some relations between these joint numerical ranges and minimal hermitian matrices.
\begin{theorem}\label{teo: intersecc mS no vacia igual a intersec JNR no nula}
	Under the previous assumptions and notations of Theorem \ref{teo mV es JNR y suma 1}, given $\V,\W$ non trivial subspaces of $\C^n$, the following statements are equivalent
	\begin{enumerate}
		\item $\mV\cap m_\W\neq \emptyset$,
		\item \label{item 2 teo relac joint num range}$	W(P_S E_1P_S,   \dots,P_S E_nP_S )\cap  W(P_\W E_1P_\W,  \dots,P_\W E_nP_\W )\neq \{0\}$
		
		(there exists a non-zero vector in the intersection), and
		\item 
		The matrix  $\lambda(P_\V-P_\W) +R$ is a minimal hermitian matrix for $R\in M^h_n(\C)$ such that $R P_\V=R P_\W=0$ and $\|R\|<\lambda$.
	\end{enumerate}
\end{theorem}

\begin{proof}
	(1) $\Rightarrow$ (2). This is apparent, since $m_\V\subset W(P_\V E_1P_\V,   \dots,P_\V E_nP_\V )$ and $m_\W\subset\allowbreak W(P_\W E_1P_\W,  \dots,P_\W E_nP_\W )$ (see \eqref{eq: mV es JNR interseccion suma 1}) and these moments are not empty and do not have the null vector.

	(2) $\Rightarrow$ (1).  If assumption  (\ref{item 2 teo relac joint num range}) holds then there exist $\rho, \mu\in\DD$ such that
	\[
	\left(\tr(P_\V E_1 P_\V\rho), \dots,\tr(P_\V E_n P_\V\rho)\right)=
	\left(\tr(P_\W E_1 P_\W\mu), \dots,\tr(P_\W E_n P_\W\mu)\right)\neq 0
	\]
	and then
	\begin{equation}
	\label{eq: igualdad de elementos del joint num range}
	\left(\tr(E_1 P_\V\rho P_\V), \dots,\tr(E_n P_\V\rho P_\V )\right)=
	\left(\tr(E_1 P_\W\mu P_\W), \dots,\tr(E_n P_\W\mu P_\W)\right).
	\end{equation}
	Then since $P_\V\rho P_\V\geq 0$ and $P_\W\mu P_\W\geq 0$ and non zero
	\begin{equation}\label{eq: igualdad de trazas de los 2 joint num range}
	\tr(P_\V\rho P_\V)=\sum_{i=1}^n \tr(E_i P_\V\rho P_\V)=\sum_{i=1}^n \tr(E_i P_\V\mu P_\V)=\tr(P_\W\mu P_\W)\neq 0
	\end{equation}
	holds.
	
	Now define
	\[
	Y=\frac{1}{\tr(P_\V \rho P_\V)}P_\V \rho P_\V\ \ \text{ and } \ \ X=\frac{1}{\tr(P_\W \mu P_\W)}P_\W \mu P_\W.
	\]
	It is apparent that $Y\geq 0$, $P_\V Y=Y$, $\tr(Y)=1$ and therefore $Y\in\DV$. Similarly $X\in\mathcal{D}_\W$.	
	
	Then considering \eqref{eq: igualdad de elementos del joint num range} and \eqref{eq: igualdad de trazas de los 2 joint num range} follows that $\Phi(Y)=\Phi(X)$ which in turn implies that $m_\V\cap m_\W\neq \emptyset$.
	
	(3) $\Leftrightarrow$ (1) This has been already proved elsewhere, see Remark \ref{rem: intersecc momentos entonces construcc mat minimal} for details.
\end{proof}\
Note that if $\V$ and $\W$ are orthogonal and satisfy any of the equivalent statements of the previous theorem then they form a support, the main object of study in \cite{soportes} (see for example Theorem 3 in that work).
\begin{corollary}
	Let $M$ be a hermitian matrix with $\|M\|=\lambda_{max}=-\lambda_{min}$, $\V=Eig_{\|M\|}$ and $\W=Eig_{-\|M\|}$ the eigenspaces of the eigenvalues $\pm\|M\|$. Then the following statements are equivalent
	\begin{enumerate}
		\item \label{a) } $(\V,\W)$ is a support
		\item \label{b) }$m_\V\cap m_\W\neq \emptyset$	 
		\item $M$ is a minimal matrix
		\item there exists a non-zero vector in the set
		$W(P_\V E_1P_\V,   \dots,P_\V E_nP_\V )\cap  W(P_\W E_1P_\W,  \dots,P_\W E_nP_\W )$.
	\end{enumerate}

\end{corollary}
\begin{proof}
	The equivalence follows after considering the definition of a support (see \cite{soportes}), Theorem \ref{teo: intersecc mS no vacia igual a intersec JNR no nula} and the equivalences between (2) and (4) in Theorem 3 of \cite{soportes}.
\end{proof}

Next we will write just $W$ and $W_\text{class}$ to denote the sets defined in \eqref{def JNR} and \eqref{def JNR clasico} respectively.

The previous results in this section related to the joint numerical ranges $W(P_S E_1 P_S, \dots,P_S E_n P_S ) \subset \left(\R_{\geq 0}\right)^n$ show that they can be used to detect minimal hermitian matrices in a very similar way we used the moment set $m_S$. 
	Nevertheless, the precise description of joint numerical ranges is not an easy task even in low dimensions. 
	
	For example, the classical joint numerical range $W_\text{class}$ of \eqref{def JNR clasico} of any triple of hermitian $n \times n$ matrices is convex if $n > 2$. Here we are interested in $n$ matrices of $n\times n$ then this property leave us only with the case of $3\times 3$ matrices. The equivalent situation for $S\in\C^3$ where the convex hull in Definition \ref{def: momento de un subespacio} is not required was settled in \cite{ammrv}.
	
	In general, the convexity of $W_{\text{class}}$ is an open problem for $n$-tuples of matrices when $n>3$ (see \cite{gutkin-jonckheere-karow, li-poon}). Similarly, in the case of subspaces $S$ with dim$(S)>3$ there exist non convex sets $\{|v|^2\in\R^n_{\geq 0}: v\in S \wedge \|v\|=1\}$ that even have not empty interior (see Remark 4 in \cite{alrv} and \cite{soportes}).
	
	Since the set $W_{\text{class}}$ may be significantly simpler than $W$, any positive result in these direction gives an easier characterization of $m_S$ (using Corollary \ref{coro Wclass convexo y mS}).
	
	
	In \cite{szymanski-weis-zyczkowski} a detailed classification of the possible joint numerical ranges of 3 matrices of $3\times 3$ is developed. In this case $W$ is a three-dimensional oval in which every one dimensional face is a segment and every two dimensional face is a filled ellipse. A characterization of only ten configurations of these segments and ellipses are possible.

We include here some general and particular properties of joint numerical ranges found in the bibliography that provide information about $m_S$.

\begin{remark}\label{rem propiedades mS a partir de W}
	Let $\{0\}\subsetneq S\subsetneq \C^n$ be a subspace of $\C^n$, $P_S$ its orthogonal projection and $E_i=e_i\cdot ( e_i )^t=e_i\otimes e_i\in M_n^h(\C)$, with $e_i$ for $i=1,\dots,n$ the standard vectors of a fixed basis $E$ as defined before.
	
	 Then, using the bibliography in joint numerical ranges, the following properties regarding $W$ (see \eqref{def JNR}), its classical version $W_\text{class}$ (see \eqref{def JNR clasico}), $W_{S,E}=W(P_S E_1 P_S,   \dots,P_S E_nP_S )$ and $m_{S,E}$ (see \eqref{def: momento de un subespacio}) hold.
	\begin{enumerate}
		\item[a) ]
		In \cite{chien-nakazato} the authors presented an example where the analogous of the Kippenhahn boundary generating curve of the numerical range $W(A_1 + i A_2 )$ for $A_1, A_2\in M_n^h$ does not hold in the case of $W$ of three $3\times 3$ hermitian matrices. 
		They consider an algebraic variety in the projective space $CP^m$ and obtain a boundary generating hypersurface that, under some hypothesis, allow the detection of elements of $W$ for any $m$ and $n$ (see Theorem 2.6 in that work). 

		\item[b) ]
		A complete generalization of the so called boundary generating curve of Kippenhahn for joint numerical ranges of hermitian matrices of any size was developed in \cite{plaumann-sinn-weis}. In that work it is proved that $W$ is the convex hull of a semi-algebraic set $T^{\sim}\subset W$ that contains the exposed points of $W$ and hence its extreme points. This subset $T^\sim$ is the euclidean closure of the union of the dual varieties of the regular points of the irreducible components of a hypersurface related to the zeroes of $\det(x_0 I + x_1 P_S E_1 P_S + \dots + x_n P_S E_n P_S ) \in \R[x_0 , x_1 , \dots , x_n ]$. Moreover, it is proved that $T^\sim$ contains the Zariski closure of the set of extreme points of $W$ (see Remarks 1.3 4) and 4.15 of \cite{plaumann-sinn-weis}).	
		Nevertheless, $T^\sim$ is not necessarily contained in the boundary of $W$.

		\item[c) ] 
		{\label{item de Remark props mS y W semialgebraico}
		 $W_{S,E}$ and $m_{S,E}$ are semi-algebraic sets of $\R^n$. 
		This follows using Theorem 1.2 and Remark 1.3 3) of \cite{plaumann-sinn-weis} and Theorem \ref{teo mV es JNR y suma 1}.
	}

		\item[d) ] 
		The study of the problem of finding all the unitary matrices $U$ such that que $m_{U(S)}=m_S$ (see Example \ref{ejemplo 1}) or of the subspaces $\V, \W$ such that $m_\V=m_\W$ is closely related to the same problems in terms of the corresponding joint numerical ranges $W_{\V,E}$ and $W_{\W,E}$. 
		
		\item[e) ]		
		The work \cite{krupnik spitkovsky} describes which matrices can produce a given $W_\text{class}(A_1, \dots, A_m)$ for $A_i\in M^h_n$, $i=1,\dots,m$. 
		More precisely, under what conditions on $m, n$ and/or the
		shape of $W_\text{class}(A_1, \dots, A_m)$ the $m$-tuple $A_1, \dots , A_m$ can be restored from it up to unitary similarity. The cases considered were $m = 2$, $n > 2$ and $n = m = 3$.

		\item[f) ] $W_\text{class}(A_1,A_2)$ is always convex and so it is $W_\text{class}(A_1,A_2,A_3)$ for $n \geq 3$ (see \cite{gutkin-jonckheere-karow, li-poon}).
		Nevertheless, the convexity of $W_\text{class}$ seems to be unanswered in the general theory of joint numerical ranges at least for big $m$.
		\item[g) ] There exist some standard numerical algorithms to generate or approximate the boundary of $W$. See for example the one detailed in page 6 of \cite{szymanski_tesis_master}.

		\end{enumerate}
	\end{remark}

\vspace{1cm}

%

	As seen in the item \eqref{item 1 remark 15} of Remark \ref{rem props nuestros jnr}, given a generic subspace $S$, the set $W_{S,E}$ can be easily described in terms of a joint numerical range of projections of rank one or as in  and \eqref{item 3 remark 15} of the same remark in terms of the coordinate subspace $\{e_i\}_{i=1}^r$ and matrices of rank one that are zero outside the $r\times r$ first block. 
	The authors did not find any reference to studies of the joint numerical range of these particular $n$-tuple of $n\times n$ matrices.


\begin{thebibliography}{XXXXXX}


\bibitem{ammrv}
Andruchow, Esteban; Mata-Lorenzo, Luis E.; Mendoza, Alberto; Recht, L\'azaro; Varela, Alejandro. {
	Minimal matrices and the corresponding minimal curves on flag manifolds in low dimension}. Linear Algebra Appl. 430 (2009), no. 8-9, 1906--1928. 
\bibitem{duranmatarecht} Dur\'an, C.E., Mata-Lorenzo, L.E. and Recht, L., { Metric geometry in homogeneous spaces of the unitary group of a C$^*$-algebra: Part I--minimal curves}, Adv. Math. 184 no. 2 (2004), 342-366.

\bibitem{plaumann-sinn-weis} Plaumann, D., Sinn, R., Weis, S. Kippenhahn's Theorem for joint numerical ranges and quantum states. SIAM Journal on Applied Algebra and Geometry, vol. 5 no 1 (2021), p. 86--113.

\bibitem{soportes} Mendoza, A., Recht, L., Varela, A. (2020). Supports for minimal hermitian matrices. Linear Algebra and its Applications, 584, 458--482. 

\bibitem{alrv}
Andruchow, Esteban; Larotonda, Gabriel; Recht, L\'azaro; Varela, Alejandro. {
	A characterization of minimal Hermitian matrices}. Linear Algebra Appl., 436 (2012) 2366--2374. 

\bibitem{bonsall-duncan} Bonsall, F. F. and Duncan, J.  {
	Numerical Ranges of Operators on Normed Spaces and of Elements of Normed Algebras}, Cambridge University Press, London, 1971.
\bibitem{weis2017} Weis, S. {On a theorem by Kippenhahn}. Preprint, arXiv:1705.00935v1 [math.AG].

\bibitem{muller} Müller, V. The joint essential numerical range, compact perturbations, and the Olsen problem, Studia Mathematica, 197  (2010), 275--290.

\bibitem{gutkin-jonckheere-karow} Gutkin, E., Jonckheere,  E. A., Karow,  M.  Convexity of the joint numerical range: topological and differential geometric viewpoints, Linear Algebra and its Applications 376 (2004), 143-171.

\bibitem{li-poon} Li, C.-K.  and Poon, Y.-T. Convexity of the joint numerical range, SIAM Journal on Matrix Analysis and Applications, 21 (2000), 668--678.

\bibitem{szymanski-weis-zyczkowski} Szyma\'nski, Konrad, Stephan Weis, and Karol {\.Z}yczkowski. Classification of joint numerical ranges of three hermitian matrices of size three. Linear algebra and its applications, 545 (2018), 148--173.


\bibitem{chien-nakazato} Chien, M. T. and Nakazato, H. Joint numerical range and its generating hypersurface. Linear algebra and its applications, 432 no 1 (2010), 173-179.

\bibitem{krupnik spitkovsky} Krupnik, N., and Spitkovsky, I.M. {
	Sets of matrices with given joint numerical range.} Linear algebra and its applications vol. 419, no 2-3 (2006), p. 569-585.







\bibitem{keeler-rodman-spitkovsky} Keeler, D. S., Rodman, L., Spitkovsky, I. M.  The numerical range of 3$\times$3 matrices. Linear algebra and its applications, 252 (1-3) (1997), 115-139.

 



\bibitem{szymanski_tesis_master} Szyma\'nski, K. Uncertainty relations and joint numerical ranges, MSc. thesis, preprint (2017) arXiv:1707.03464v1 [quant-ph].



%



\end{thebibliography}
\end{document}